\documentclass[11pt]{article}

\usepackage{amsmath}
\usepackage{amssymb}
\usepackage{eucal}

\begin{document}

\baselineskip 16pt

\title{Finite groups whose  $n$-maximal subgroups are   $\sigma$-subnormal \thanks{Research is supported by
a NNSF grant of China (Grant \# 11371335) and Wu Wen-Tsun Key Laboratory of Mathematics of Chinese Academy of Sciences.} }

\author{Wenbin Guo\\
{\small Department of Mathematics, University of Science and
Technology of China,}\\ {\small Hefei 230026, P. R. China}\\
{\small E-mail:
wbguo@ustc.edu.cn}\\ \\
{ Alexander  N. Skiba}\\
{\small Department of Mathematics,  Francisk Skorina Gomel State University,}\\
{\small Gomel 246019, Belarus}\\
{\small E-mail: alexander.skiba49@gmail.com}}

\date{}
\maketitle

\begin{abstract} Let $\sigma =\{\sigma_{i} | i\in I\}$ be some partition of the set of all primes $\Bbb{P}$.
A set  ${\cal H}$  of subgroups of $G$ is said to be a  \emph{complete Hall $\sigma $-set} of $G$   if   every  member $\ne 1$ of
 ${\cal H}$ is a Hall $\sigma _{i}$-subgroup  of $G$, for some $i\in I$, and $\cal H$ contains exact one Hall  $\sigma _{i}$-subgroup of $G$
 for every  $\sigma _{i}\in  \sigma (G)$. A  subgroup  $H$ of $G$ is said to be:   \emph{$\sigma$-permutable} or \emph{$\sigma$-quasinormal} in $G$ if $G$ possesses a complete Hall $\sigma$-set set ${\cal H}$ such that $HA^{x}=A^{x}H$ for all  $A\in {\cal H}$ and  $x\in G$:   \emph{${\sigma}$-subnormal} in $G$   if there is a subgroup chain  $A=A_{0} \leq A_{1} \leq \cdots \leq A_{t}=G$  such that  either $A_{i-1}\trianglelefteq A_{i}$ or $A_{i}/(A_{i-1})_{A_{i}}$ is    a finite $\sigma_{i}$-group for some $\sigma_{i}\in \sigma$   for all $i=1, \ldots t$.

If $M_n < M_{n-1} < \ldots < M_1 < M_{0}=G$, where $M_i$ is a maximal subgroup of  $M_{i-1}$, $i=1,2,\ldots ,n$, then $M$ is said to be an \emph{$n$-maximal subgroup} of $G$.
   If each $n$-maximal subgroup  of $G$ is $\sigma$-subnormal ($\sigma$-quasinormal, respectively) in $G$ but, in the case $ n > 1$,  some $(n-1)$-maximal subgroup is not $\sigma$-subnormal (not $\sigma$-quasinormal, respectively)) in $G$, we   write $m_{\sigma}(G)=n$ ($m_{\sigma q}(G)=n$, respectively).

In this paper, we show that the parameters  $m_{\sigma}(G)$ and $m_{\sigma q}(G)$  make possible to bound the $\sigma$-nilpotent length $ \
l_{\sigma}(G)$ (see below the definitions of the terms employed),  the rank $r(G)$  and  the number $|\pi (G)|$ of all distinct primes dividing the order  $|G|$ of a finite soluble  group $G$.
  We also  give conditions under which a finite group  is  $\sigma$-soluble or $\sigma$-nilpotent, and
   describe the structure of a finite soluble group $G$  in the case when $m_{\sigma}(G)=|\pi (G)|$.
   Some known results are generalized.

\end{abstract}


\let\thefootnote\thefootnoteorig

\section{Introduction}

Throughout this paper, all groups are finite and $G$ always denotes a finite group.  Moreover,  $\mathbb{P}$ is the set of all primes, $\pi \subseteq  \Bbb{P}$ and  $\pi' =  \Bbb{P} \setminus \pi$. If $n$ is an integer, the symbol $\pi (n)$ denotes the set of all primes dividing $|n|$; as usual,  $\pi (G)=\pi (|G|)$, the set of all  primes dividing the order of $G$.
Let   $A$ and $B$ be  subgroups of $G$. We say that $A$ forms an \emph{irreducible pair  with $B$} if $AB=BA$ and $A$ is a maximal subgroup of $AB$.

In what follows, $\sigma =\{\sigma_{i} | i\in I\subseteq \mathbb{N} \}$ is  some partition of $\mathbb{P}$, that is, $\Bbb{P}=\cup_{i\in I} \sigma_{i}$ and $\sigma_{i}\cap
\sigma_{j}= \emptyset  $ for all $i\ne j$.    By analogy with the notations $\pi (n)$ and $\pi (G)$, we put $\sigma (n) =\{\sigma_{i} |\sigma_{i}\cap \pi (n)\ne  \emptyset  \}$  and  $\sigma (G)=\sigma (|G|)$.

In the mathematical practice, we often deal with the following two special partitions of  $\Bbb{P}$:
 $\sigma =\{\{2\}, \{3\},  \ldots \}$ and    $\sigma =\{\pi, \pi'\}$.

A group $G$ is called:  \emph{$\sigma$-primary} \cite{1} if  $|\sigma (G)|\leq 1$, that is, $G$ is a $\sigma _{i}$-group for some $i$;
 \emph{$\sigma$-nilpotent} or \emph{$\sigma$-decomposable} (Shemetkov \cite{bookShem}) if $G=A_{1}\times \cdots \times A_{r}$ for some
 $\sigma$-primary groups $A_{1}, \ldots , A_{r}$;  \emph{$\sigma$-soluble} \cite{1} if every chief factor of $G$ is  $\sigma$-primary.
 By the \emph{$\sigma$-nilpotent length} (denoted by $l_{\sigma}(G)$)  of a  $\sigma$-soluble group $G$ we mean the length of the shortest normal chain of $G$ with $\sigma$-nilpotent factors.

  The group   $G$ is nilpotent if and only if  $G$ is $\sigma$-nilpotent where  $\sigma =\{\{2\}, \{3\},  \ldots \}$; $G$ is $\pi$-decomposable, that is,
 $G=O_{\pi}(G) \times  O_{\pi'}(G)$ if and only if $G$ is $\sigma$-nilpotent where  $\sigma =\{\pi, \pi'\}$.

 The $\sigma$-nilpotent groups have many applications in the formation theory  \cite{bookShem, 19, 20, Bal-Ez} (see also the recent papers \cite{1, wiel} and the survey \cite{30}), and such groups  are exactly the  groups in which every subgroup is $\sigma$-subnormal in the sense of the following definition (see Proposition 3.4 below).

{\bf Definition 1.1.} A  subgroup $A$ of $G$ is called \emph{${\sigma}$-subnormal} in $G$ \cite{1}  if there is a subgroup chain  $$A=A_{0} \leq A_{1} \leq \cdots \leq
A_{t}=G$$  such that  either $A_{i-1}\trianglelefteq A_{i}$ or $A_{i}/(A_{i-1})_{A_{i}}$ is  ${\sigma}$-primary for all $i=1, \ldots , t$.

It should be noted that  the  idea of this concept originates from the paper of Kegel \cite{5} where it was discussed one another generalization of nilpotency.

There are some motivations for introducing and study of $\sigma$-subnormal subgroups. First of all, the set of all $\sigma$-subnormal subgroups of $G$ forms a sublattice of the lattice of all subgroups of $G$ (see Proposition 2.5 below). This fact is a generalizations of the classical Wielandt's result which states that
the set of all subnormal subgroups of $G$ forms a sublattice of the lattice of all subgroups of $G$. Thus, the $\sigma$-subnormal subgroups are very convenient for applications. The first among such applications were obtained in \cite{1, 2}.

Recall that a set  ${\cal H}$
 of subgroups of $G$ is said to be a  \emph{ complete
Hall $\sigma $-set} of $G$   if   every  member $\ne 1$ of
 ${\cal H}$ is a Hall $\sigma _{i}$-subgroup
 of $G$, for some $i\in I$, and $\cal
 H$ contains exact one Hall  $\sigma _{i}$-subgroup of $G$
 for every  $\sigma _{i}\in  \sigma (G)$. A  subgroup  $H$ of $G$ is said to be   \emph{$\sigma$-permutable} or \emph{$\sigma$-quasinormal} in $G$ \cite{1} if $G$ possesses a complete Hall $\sigma$-set
set ${\cal H}$ such that $HA^{x}=A^{x}H$ for all  $A\in {\cal H}$ and  $x\in G$. In particular, $H$ is said to be   \emph{$S$-permutable} or \emph{$S$-quasinormal} in $G$ \cite{prod, GuoII} if $H$ permutes with all Sylow subgroups of $G$.

If $H$ is an $S$-quasinormal subgroup of $G$, then $H$ is subnormal in $G$ (Kegel) and $H/H_{G} $ is nilpotent (Deskins) \cite[Theorem  1.2.14]{prod}). But in general, if $H$ is $\sigma$-quasinormal in  $G$, then $H$ is not necessary subnormal and $H/H_{G}$ may be non-nilpotent \cite{1}.
Nevertheless, in this case when $H$ is $\sigma$-quasinormal,  $H$ is $\sigma$-subnormal in  $G$ and $H/H_{G}$ is $\sigma$-nilpotent \cite[Theorem B]{1}.

Examples and some other applications of   $\sigma$-subnormal subgroups and $\sigma$-quasinormal subgroups were discussed in \cite{1, 30, 2}. In this paper, we consider some applications of such  subgroups in the theory of $n$-maximal subgroups.

  Recall that if $M_n < M_{n-1} < \ldots < M_1 < M_{0}=G$ (*), where $M_i$ is a maximal subgroup of  $M_{i-1}$ for all $i=1, \ldots ,n$, then the chain   (*)  is said to be  a \emph{maximal chain of $G$ of length $n$} and $M_n $ ($n > 0$),  is an \emph{$n$-maximal subgroup}  of $G$.

   If each $n$-maximal subgroup  of $G$ is $\sigma$-subnormal ($\sigma$-quasinormal, respectively) in $G$ but, in the case $ n > 1$,  some $(n-1)$-maximal subgroup is not $\sigma$-subnormal (not $\sigma$-quasinormal, respectively)) in $G$, we   write $m_{\sigma}(G)=n$ ($m_{\sigma q}(G)=n$, respectively).

Note that $m_{\sigma}(G)=1=m_{\sigma q}(G)$ if and only if $G$ is $\sigma$-nilpotent by Proposition 3.4 below. We also show (see Corollaries 1.7 and 1.8 below)  that  $m_{\sigma}(G)=2$ if and only if $G$ is a  Schmidt group $G$ with abelian Sylow subgroups such that     $|\sigma (G)|=|\pi (G)|$,  and  $m_{\sigma q}(G)=2$ if and only if $G$ is a supersoluble group with $m_{\sigma}(G)=2$.  Finally, note that every group with $m_{\sigma}(G)=3$
is $\sigma$-soluble (see Theorem 1.4 below), and there are non-soluble groups, for example the  alternating group $A_{5}$ of degree 5,    with $m_{\sigma q}(G)=4$.

If $G$ is soluble, the parameters  $m_{\sigma}(G)$ and $m_{\sigma q}(G)$  make possible to bound the $\sigma$-nilpotent length $ l_{\sigma}(G)$,  the rank $r(G)$ and  the number $|\pi (G)|$ of all distinct primes dividing $|G|$.

Recall  that the \emph{rank} $r(G)$ of a soluble group $G$ is the maximal integer $k$ such that $G$ has a chief factor of order $p^{k}$ for some prime $p$ (see \cite[p. 685]{hupp}).

{\bf Theorem  1.2.} {\sl Suppose that $G$ is $\sigma$-soluble and let $\cal H$
be   a complete Hall $\sigma$-set of $G$.  Then the following statements hold.}

(i) {\sl   If $G$ is soluble but it is not $\sigma$-nilpotent and
 $r(H) \leq r\in \Bbb{N}$ for all $H\in \cal H$, then
 $r(G) \leq m_{\sigma q}(G) + r-2$.}

(ii) {\sl   $l_{\sigma}(G)\leq m_{\sigma}(G)$.}

(iii) {\sl  If $G$ is soluble but it is not $\sigma$-nilpotent, then $|\pi (G)|\leq m_{\sigma}(G)$.}

Now, let's consider some applications of Theorems 1.2.

The relationship between $n$-maximal subgroups (where  $n>1$) of a group $G$ and the structure of $G$ was studied by many authors (see, in particular, the recent  papers \cite{li1, guo1, li2, KM1, KM2, KS3, Kovaleva1, K2} and Chapter 4 in the book \cite{GuoII}). One of the first results in  this direction were obtained by Huppert \cite{HupI}. In fact,
Huppert proved \cite{HupI} that: if every 2-maximal
 subgroup of $G$ is normal in $G$, then $G$ is supersoluble;    if
 every 3-maximal
 subgroup of $G$ is normal in $G$, then $G$ is soluble of rank $r(G)$
 at most two.  The first of these two results  was generalized by Agrawal
\cite{Agr} :  {\sl If every 2-maximal
 subgroup $L$ of $G$ is $S$-quasinormal in $G$,  then $G$ is supersoluble}.
In the universe of all soluble groups the  both Huppert's observations  and
some similar results in  \cite{{yanko}}
 are special cases of the following general result  (Mann \cite{mann}):
 {\sl If $G$ is soluble and every
 $n$-maximal subgroup $L$ of $G$ $(n > 1)$ is quasinormal (that is, $L$  permutes with
 all  subgroups  of $G$), then  $r(G) \leq n-1$.}

 In the case  $\sigma =\{\{2\}, \{3\}, \ldots \}$ we get from Theorem 1.2(i) the following generalization of the last  of these results.

{\bf Corollary  1.3.}  {\sl Suppose that $G$ is    soluble
 and each $n$-maximal subgroup of $G$ $(n > 1)$ is $S$-quasinormal in $G$. Then
 $r(G) \leq n-1$.}

 The following theorem allows us to obtain the above mentioned result of Agrawal.

{\bf Theorem  1.4}.  (i) {\sl If in every maximal
 chain $M_{3} < M_{2} < M_{1} <  M_{0}=G $ of $G$ of  length $3$,  one of
  $M_{3}$, $M_{2}$ and $ M_{1}$ is
$\sigma$-subnormal in $G$, then  $G$ is $\sigma$-soluble.  }

(ii) {\sl If $1  < m_{\sigma }(G)\leq 3$, then $G$ is soluble.}

{\bf Corolary  1.5} (Spencer \cite{spencer}).  {\sl If in  every maximal
 chain $M_{3} < M_{2} < M_{1} <  M_{0}=G $ of $G$ of  length $3$, one of
  $M_{3}$,  $M_{2}$  and $ M_{1}$ is
subnormal in $G$, then  $G$ is  soluble.  }

{\bf Corollary 1.6} (Mann \cite{mann}). {\sl If   every $3$-maximal
  subgroup of $G$ is
subnormal in $G$, then  $G$ is  soluble.  }

Recall that $G$ is called a \emph{Schmidt group} if $G$ is not nilpotent but every proper subgroup of $G$ is nilpotent.

{\bf Corollary 1.7.} {\sl The equality   $m_{\sigma}(G)=2$ is true  if and only if $G$ is a  Schmidt group $G$ with abelian Sylow subgroups such that     $|\sigma (G)|=|\pi (G)|$. }

{\bf Corollary 1.8.} {\sl The equality   $m_{\sigma q}(G)=2$ is true  if and only if $G$ is a supersoluble group with $m_{\sigma}(G)=2$. }

From Corollaries 1.8 and Theorem 1.4, we get

{\bf Corollary 1.9} (see Agrawal \cite{Agr}  or
  Theorem 6.5  in \cite[Ch.1]{We})). {\sl If   every $2$-maximal
  subgroup of $G$ is
$S$-quasinormal in $G$, then  $G$ is  supersoluble.  Moreover, if $|\pi (G)|> 2$, then $G$ is nilpotent.}

From Theorem 1.2(iii) we know that for every soluble but non-$\sigma$-nilpotent  group $G$ we have $|\pi (G)|\leq m_{\sigma}(G)$. In the case when  $|\pi (G)|= m_{\sigma}(G)$
the structure of such a group  $G$ can be described  completely as follows.

{\bf Theorem  1.10.} {\sl Suppose that $G$ is  soluble.  Then $m_{\sigma}(G)=|\pi (G)|$ if and only if $G$ is a group of one of the following two types:}

(i) {\sl $G$ is a $p$-group for some prime $p$.}

(ii) {\sl  $G=D\rtimes M $, where  $D=G^{{\mathfrak{N}}_{\sigma}}$ is an abelian Hall subgroup of $G$, and the following hold:}

(a) {\sl  Every non-$\sigma$-subnormal Sylow subgroup $P_{1}$ of  $G$
   is cyclic and  the maximal subgroup of $P_{1}$ is $\sigma$-subnormal in $G$. Moreover, if $P_{1},  \ldots , P_{n}$ is  a  Sylow  basis of $G$, then $P_{2},  \ldots , P_{n}$   are elementary abelian and   and $P_{1}$  forms an irreducible pairs with all such subgroups;   if
 $\{H_{1}, \ldots , H_{t}\}$ is a complete Hall $\sigma$-set of $G$,   $P_{1}\leq H_{1}$ and $P_{1}$ is not of prime order, then $H_{2},  \ldots , H_{t}$ are normal in $G$.}

(b) {\sl Some Sylow subgroup of $M$ is not $\sigma$-subnormal in $G$. Hence $M$ acts irreducibly on every Sylow subgroup of $D$.}

(c) {\sl If $G$ possesses at least two non-$\sigma$-subnormal non-isomorphic Sylow subgroups, then all non-$\sigma$-subnormal Sylow subgroups are of prime order. }

(d) {\sl If $P$ is a non-$\sigma$-subnormal Sylow subgroup of $G$,  $P$ is a  $\sigma _{i}$-group and $V$ is the maximal subgroup  of $P$, then $|G:N_{G}(V)|$  is a $\sigma _{i}$-number. }

In this theorem $G^{{\mathfrak{N}}_{\sigma}}$ denotes the $\sigma$-nilpotent residual of $G$, that is, the intersection of all normal subgroups $N$ of $G$ with $\sigma$-nilpotent quotient $G/N$.

In the case when $\sigma =\{\{2\}, \{3\}, \ldots \}$,  from Theorems  1.2 and 1.10 we get the following known result.

{\bf Corollary 1.11} (Mann \cite{mann}). {\sl Suppose that $G$ is a   soluble group and each $n$-maximal subgroup of $G$ is subnormal.
If  $n \leq   |\pi (G)|$, then $G$ is either of the following type:}

(a) {\sl $G$ is nilpotent.}

(b) {\sl $G=HN$, where,}

(i) {\sl $N$ is a normal abelian Hall subgroup, and all Sylow subgroups of $N$ are elementary abelian. }

(ii) {\sl $H$ is a cyclic Hall subgroup, and $|H|$ is either a prime power or square-free number.}

(iii)  $(|N|, |H|)=1.$

(iv) {\sl  If $H_{p}$ is a Sylow subgroup of $H$ and  $N_{q}$ is a Sylow subgroup of $N$, then $H_{p}$ induces in $N_{q}$ an irreducible automorphism group of order $p$ or 1. In the latter case, $|N_{q}|$ =q.}

{\sl Conversely, a group of type (a) or (b) has each $n$-maximal subgroup subnormal.}

{\bf Corollary 1.12} (Mann \cite{mann}). {\sl Suppose that $G$ is a   soluble group and each $n$-maximal subgroup of $G$ is subnormal.
If  $n <  |\pi (G)|$, then $G$ is nilpotent.}

We prove Theorems 1.2, 1.4 and 1.10 in Sections 4, 5 and 6, respectively. But before them, as  preparatory steps,  we prove in Section 2 that the set of
all $\sigma$-subnormal subgroups of $G$ forms a sublattice of the lattice of all subgroups of $G$,  and we collect  in Section 3 some needed
  properties of    $\sigma$-soluble and $\sigma$-nilpotent groups.

All unexplained notation and terminology are standard. The reader is referred to \cite{DH},  \cite{Bal-Ez} and \cite{GuoII} if  necessary.

\section{The lattice ${\cal L}_{\sigma}(G)$ of all $\sigma$-subnormal subgroups}

It is not difficult to show that the intersection of any two $\sigma$-subnormal subgroups of $G$ is also $\sigma$-subnormal in $G$ (see Lemma 2.2(1)(2) below). It is well-known that any partially ordered set with 1 in which  there is the greatest lower bound for each its non-empty subset is a lattice.
 Hence the set ${\cal L}_{\sigma}(G)$ of all $\sigma$-subnormal subgroups of $G$ is  a lattice. In this section, we show that ${\cal L}_{\sigma}(G)$ is a sublattice of the lattice of all subgroups of $G$.

We use ${\mathfrak{S}}_{\sigma}$   and ${\mathfrak{N}}_{\sigma}$  to denote the class  of all $\sigma$-soluble groups and the  class  of all $\sigma$-nilpotent groups, respectively.

{\bf Lemma 2.1 } (See Lemma 2.5 in \cite{1}).
 {\sl  The class  ${\mathfrak{N}}_{\sigma}$ is closed under taking direct
products, homomorphic images and  subgroups. Moreover, if  $H$ is a normal
subgroup of $G$ and  $H/H\cap \Phi (G)$ is $\sigma$-nilpotent, then
$H$ is $\sigma$-nilpotent.    }

In what follows,  $\Pi$ is always supposed to be a non-empty   subset of the set $\sigma$  and $\Pi'=\sigma\setminus \Pi$.
 We say that: a natural number $n$ is  a \emph{${\Pi}$-number} if $\sigma (n)\subseteq
{\Pi}$; $G$ is a  \emph{${\Pi}$-group} if $|G|$ is a ${\Pi}$-number;  $G$ is \emph{$\sigma $-perfect}
 if $G^{{\frak{N}}_{\sigma}}=G$.

We call the product of all normal $\sigma$-nilpotent  subgroups of $G$ the  \emph{$\sigma$-Fitting subgroup}  of $G$ and denote it by $F_{\sigma}(G)$.

{\bf Lemma  2.2.} {\sl Let  $A$,  $K$ and $N$ be subgroups of  $G$.
 Suppose that   $A$
is $\sigma $-subnormal  in $G$ and $N$ is normal in $G$.
    }

(1) {\sl $A\cap K$    is  $\sigma $-subnormal in
$K$}.

(2) {\sl If  $K$ is  a $\sigma $-subnormal subgroup of   $A$,
then $K$ is $\sigma $-subnormal   in $G$}.

(3) {\sl If $N\leq K$ and $K/N$ is $\sigma $-subnormal  in $G/N$, then $K$
is $\sigma $-subnormal  in $G$}.

(4) {\sl If $K\leq A$ and $A$ is $\sigma $-nilpotent, then $K$
is $\sigma$-subnormal    in $G$.}

(5) {\sl If $N$ is a $\sigma _{i}$-group, for some $i$,
then $N\leq N_{G}(O^{\sigma _{i}}(A))$.}

(6) {\sl If $K\leq E\leq G$, where $K$ is $\sigma $-subnormal  in   $E$,
 then $KN/N$ is
$\sigma $-subnormal  in $NE/N$. }

(7) {\sl If  $A$ is $\sigma $-perfect, then
 $A$ is subnormal in $G$. }

(8)   {\sl Suppose that $N$ is   the product
of some  minimal normal subgroups of $G$ and    $N$ is  not $\sigma $-primary.
  Suppose also that    $G=AN$, $N$ is non-abelian and
 all composition factors of $N$ are isomorphic. Then $N\leq N_{G}(A)$. }

(9) {\sl If    $A $
  is a $\sigma _{i}$-group,  then  $A\leq O_{\sigma _{i}}(G)$.
 Hence if $A$ is $\sigma$-nilpotent, then $A\leq  F_{\sigma}(G)$. }

{\bf Proof.}   Statements (1)--(5)  follow from  Lemma 2.6 in \cite{1}.

(6)   By hypothesis, there is  a chain  $K=K_{0} \leq
K_{1} \leq \cdots \leq K_{n}=E$ such that
either $K_{i-1} \trianglelefteq K_{i}$
  or $K_{i}/(K_{i-1})_{K_{i}}$ is  $\sigma $-primary  for all $i=1, \ldots , n$.
   Consider the chain $$KN/N=K_{0}N/N \leq K_{1}N/N \leq \cdots \leq
K_{n}N/N=EN/N.$$ Assume that $K_{i-1}N/N$   is not normal in $K_{i}N/N$.
Then $L=K_{i-1}$   is not normal in $T=K_{i}$ and so $T/L_{T}$
is $\sigma $-primary.
Then $$(T/L_{T})/(L_{T}(T\cap N)/L_{T})=(T/L_{T})
/((T\cap NL_{T})/L_{T})\simeq $$$$T/(T\cap NL_{T})\simeq
TN/L_{T}N\simeq (TN/N)/(L_{T}N/N)$$  is $\sigma $-primary.
 But $L_{T}N/N\leq (LN/N)_{TN/N}$.
Hence $(TN/N)/(LN/N)_{TN/N}$ is $\sigma $-primary.
This shows that $KN/N$  is  $\sigma $-subnormal in
$NE/N$.

  By hypothesis, there is a  chain  $A=A_{0} \leq
A_{1} \leq \cdots \leq A_{r}=G$ such that
either $A_{i-1} \trianglelefteq A_{i}$
  or $A_{i}/(A_{i-1})_{A_{i}}$ is  $\sigma $-primary  for all $i=1, \ldots
, r$.   Let   $M=A_{r-1}$.
  We can assume without loss of generality that $M\ne G$.

(7)    Assume that this assertion is false and let $G$ be a counterexample of minimal order.
First we show that $A\leq  M_{G}$. This is clear  if $M$ is normal in $G$. Now
assume that $G/M_{G}$ is a $\sigma _{i}$-group. Then
 from the isomorphism $AM_{G}/M_{G}\simeq A/A\cap M_{G}$ and $A=A^{\frak{N}_{\sigma}}$ we get that
  $A=O^{\sigma _{i}}(A)\leq A\cap
M_{G}$, so $A\leq M_{G}$.
The choice of $G$ implies that $A$ is subnormal in $M$, so $A$ is subnormal in $M_{G}$ by
 Assertion (1). Therefore     $A$ is subnormal in $G$.

(8) Assume that this assertion is false. First note that
 $N\nleq M$ since $G=AN$ and $M < G$ by hypothesis.

 If  $M$ is not normal in $G$, then  $G/M_{G}$ is  $\sigma $-primary and  so
  $N\simeq NM_{G}/M_{G}$ is  $\sigma $-primary, which contradicts the hypothesis.
 Hence  $M$ is normal in $G$ and so
 $N=(N\cap M)\times N_{0}$,
 where  $N\cap M$ and $ N_{0}$
 are normal in $G$. Then $N\cap M \leq \text{Soc} (M)$ by  \cite[A, 4.13(c)]{DH}
 and $M=M\cap AN=A(N\cap M)$. The choice of $G$ implies that $N\cap M\leq N_{M}(A)\leq N_{G}(A)$.
On the other hand, $N_{0}\cap M=1$ and  so $N_{0}\leq C_{G}(M) \leq C_{G}(A)$. Thus $N\leq
N_{G}(A)$, a contradiction. Hence we have (8).

(9) Assume that this assertion      is false and let $G$ be a counterexample
of  minimal order. Then $1 <  A < G$. Let $D=O_{\sigma _{i}}(G)$,  $R$ be
  a minimal normal subgroup of $G$ and $O/R=O_{\sigma_{i}}(G/R)$.   Then
the choice of $G$ and Assertion (6)  imply that   $AR/R\leq O/R$, so $A\leq O$.
 Therefore   $R\nleq D $, so $D=1$  and  $A\cap R < R$.
 Suppose that $L= A\cap  R\ne 1$.  Assertion (1)  implies that $L$
is $\sigma $-subnormal in $R$.   If $R < G$, the choice of $G$ implies that $L\leq
O_{\sigma _{i}}(R)\leq   D$ since $O_{\sigma _{i}}(R)$ is a characteristic subgroup of $R$.
   But then $D\ne 1$,
 a contradiction. Hence $R=G$ is a simple group, which is  also
 impossible since $1 <  A < G$.  Therefore $R\cap
A=1$.  If  $O < G$, the choice of $G$ implies that
$A\leq O_{\sigma _{i}}(O)\leq   D=1$. This contradiction shows that
  $G/R=O/R$ is a $\sigma _{i}$-group.
Hence $R$ is the unique minimal normal subgroup of $G$. Moreover, Lemma 2.1  implies that
 $R\nleq\Phi (G)$,  so $C_{G}(R)\leq R$ by \cite[A, 15.2]{DH}.

Now we show that $A\subseteq R$. First assume that  $R$ is   $\sigma$-primary. Then  $R$ is a ${\sigma _{j}}$-group for some
 ${\sigma _{j}}\in \sigma\setminus \{{\sigma _{i}}\}$ and so $O^{\sigma _{j}}(A)=A$.
 Therefore
$R\leq N_{G}(A)$ by Assertion (5). Consequently, $A\leq C_G(R)\leq R$.  Now assume that $R$ is not
$\sigma $-primary. Then $R$ is not
abelian.   Hence  $R$ is the product of some minimal normal subgroups of $RA$ by \cite[A, 4.13(c)]{DH}.
Hence
 $R\leq N_{RA}(A)$ by  Assertion (8). Then $AR=A\times R$ and so also
$A\leq C_{G}(R)\leq R$. This contradiction completes the proof the fact
that $A\leq O_{\sigma _{i}}(G)$. Now assume that $A$ is
$\sigma$-nilpotent. Then $A=A_{1}\times \cdots \times A_{t}$, where
$A_{1}, \ldots , A_{t}$  are $\sigma$-primary groups.
Since $A$ is $\sigma$-subnormal in $G$, every factor
$A_{i}$ is $\sigma$-subnormal in $G$. Hence $A_{i}$
 is contained in  $O_{\sigma _{j}}(G)$ for some $j=j(A_{i})$,  so
$A_{i}\leq F_{\sigma}(G)$.  Thus $A\leq F_{\sigma}(G)$.

The lemma is proved.

   Let $$F_{0\sigma}(G) \leq  F_{1\sigma}(G)   \leq \cdots \leq
 F_{i\sigma}(G) \leq \cdots $$ be
 the \emph{upper ${\sigma}$-nilpotent series}  of $G$, that is,  $F_{0\sigma}(G)=1$
and $$F_{i\sigma}(G)/F_{(i-1)\sigma}(G)=F_{\sigma}(G/F_{(i-1)\sigma}(G))$$
for $i > 0$. If $n$ is the smallest integer such that $F_{n\sigma}(G)=G$,
then  $n$  coincides with    the ${\sigma}$-nilpotent length of $G$.

We use ${\mathfrak{N}}^{n}_{\sigma}$    to denote the class  of all
$\sigma$-soluble   groups $G$ such that $l_{\sigma}(G)\leq n$ ($n > 0)$.

{\bf Lemma  2.3. }  {\sl The following holds:}

   (i) {\sl If a non-empty class ${\mathfrak{F}}$ of groups is
 closed under taking direct
products, homomorphic images and  subgroups, then the class
 ${\mathfrak{N}}_{\sigma}{\mathfrak{F}}$ is also  closed under taking direct
products, homomorphic images and  subgroups. Moreover, if
 $G/\Phi (G)\in  {\mathfrak{N}}_{\sigma}{\mathfrak{F}}$,
 then
$G\in {\mathfrak{N}}_{\sigma}{\mathfrak{F}}$.   }

 (ii) {\sl  The class ${\mathfrak{N}}^{n}_{\sigma}$   is closed under taking direct
products, homomorphic images and  subgroups. Moreover, if
  $G/\Phi (G)\in {\mathfrak{N}}^{n}_{\sigma}$, then
$G\in {\mathfrak{N}}^{n}_{\sigma}$.    }

{\bf Proof. } (i)   This assertion can be proved by the direct
calculations.   For instance, if   $G/\Phi (G)\in  {\mathfrak{N}}_{\sigma}{\mathfrak{F}}$ and
$N/\Phi (G)$ is a normal subgroup of $G/\Phi (G)$ such that  $(G/\Phi
(G))/(N/\Phi (G))\simeq G/N \in  {\mathfrak{F}}$ and  $N/\Phi (G)$ is $\sigma$-nilpotent,
 then  $N$ is $\sigma$-nilpotent by Lemma  2.1. Hence $G\in {\mathfrak{N}}_{\sigma}{\mathfrak{F}}$.

(ii)  In the case when $n=1$, the assertion follows from Lemma 2.1. Now assume
 that  $ n >1$ and  the assertion is true  for $n-1$.
 It is not difficult to show that  $${\mathfrak{N}}^{n}_{\sigma}
={\mathfrak{N}}_{\sigma} {\mathfrak{N}}^{n-1}_{\sigma}.$$ Therefore this
assertion  is a corollary of  (i).

The lemma is proved.

{\bf  Proposition  2.4.} {\sl  Let $A$ be a $\sigma$-subnormal
 subgroup of  $G$. If   $A$
 is  $\sigma$-soluble and  $l_{\sigma}(A)\leq n$, then    $A\leq F_{n\sigma}(G)$.}

{\bf  Proof.} First note that in view of Lemma  2.2(9), $F_{\sigma}(A)\leq
F_{\sigma}(G)$.  Hence
 $l_{\sigma}(AF_{\sigma}(G)/F_{\sigma}(G))= l_{\sigma}(A/ A\cap F_{\sigma}(G)) \leq n-1$
and  so  by induction we have   $AF_{\sigma}(G)/F_{\sigma}(G)\leq
F_{n-1 \sigma}(G/F_{\sigma}(G))=F_{n\sigma}(G)/F_{\sigma}(G)$.  Therefore
$A\leq F_{n\sigma}(G)$.    The lemma is proved.

{\bf Proposition 2.5.} {\sl The set of all $\sigma$-subnormal subgroups $A$ of $G$ with $l_{\sigma}(A)\leq n$ forms   a sublattice of the lattice of all subgroups of $G$. }

{\bf Proof.}  In view of Lemma 2.3, Proposition 2.4 and Statements (1) and (2)
 of Lemma 2.2,
 we need only to show that
 if $A$ and $B$ are
$\sigma $-subnormal    subgroups of $G$, then $\langle A, B\rangle$ is
 $\sigma $-subnormal  in $G$.

Assume that this
  is false and let $G$ be a counterexample of
 minimal     order.   Then $A\ne 1\ne B$  and $\langle A, B\rangle  \ne G$.
 Let $R$ be a minimal normal subgroup of $G$.

(1) {\sl $\langle A, B\rangle R=G$. Hence $\langle A, B\rangle _{G}=1$}

Suppose that $L=\langle A, B\rangle R\ne G$.
Lemma   2.2(1) implies that  $ A  $ and $B$ are
 $\sigma $-subnormal in $L$. The
choice of $G$ implies that  $\langle A, B\rangle$ is
 $\sigma $-subnormal  in $L$.  On the
other hand, $$L/R=\langle A, B\rangle R/R=\langle AR/R, BR/R\rangle,$$
where  $ AR/R  $ and $BR/R$ are
 $\sigma$-subnormal  in $G/R$ by Lemma  2.2(6),
 so the choice of $G$  implies that $L/R$ is
 $\sigma$-subnormal  in $G/R$ and so $L$ is
 $\sigma$-subnormal  in $G$ by Lemma
2.2(3). But then $\langle A, B\rangle$ is
 $\sigma $-subnormal in $G$ by Lemma
2.2(2).    This contradiction shows that we have (1).

(2) {\sl If  $S$ is a non-identity characteristic subgroup of $C$, where
$C\in \{A, B\}$, then $R\nleq N_{G}(C)$.}

Indeed, if $R\leq N_{G}(C)$, then  $R\leq N_{G}(S)$ and so
 $$S^{G}= S^{\langle A, B\rangle R}=S^{\langle A, B\rangle }\leq
C^{\langle A, B\rangle}\leq \langle A, B\rangle_{G}=1,$$ a contradiction.

(3)  {\sl   $R$  is a  $\sigma _{i}$-group   for some  $i\in I$.}

Suppose that this is false. Then $R$
is non-abelian, which implies that $R$ is the product of some minimal
normal subgroups of $RA$ by \cite[A, 4.13(c)]{DH}.   Hence $R\leq N_{RA}(A)\leq N_G(C)$ by Lemma  2.2(8), contrary to Claim (2).
Hence we have (3).

{\sl Final contradiction.} First we show that   $A$ and $B$     are
$\sigma _{i}$-groups. Indeed, since $R$ is a
   $\sigma _{i}$-group  by Claim (3), $R\leq N_{G}(O^{\sigma _{i}}(A))$ by Lemma  2.2(5).
 But   $O^{\sigma _{i}}(A)$ is a
characteristic subgroup of  $A$, so   $O^{\sigma _{i}}(A)=1$ by Claim (2). Hence
$A$ is   a  $\sigma _{i}$-group.  Similarly one can get that $B$ is a   $\sigma _{i}$-group.
Therefore $\langle A, B \rangle\leq  O_{\sigma _{i}}(G)$ by Lemma 2.2(9).
Hence $\langle A, B \rangle$ is $\sigma$-subnormal in $G$ by Lemma 2.2(4).
But this  contradicts the choice of $G$. The proposition  is proved.

{\bf   Corollary 2.6} (H. Wielandt \cite[A, 14.4]{DH}). {\sl The set of all subnormal subgroups of $G$ forms a subllatice of the lattice of all subgroups of $G$.}

\section{Some properties of $\sigma$-soluble and $\sigma$-nilpotent groups}

The direct calculations show that the following lemma is true

{\bf Lemma 3.1.}  {\sl   The class  ${\mathfrak{S}}_{\sigma}$ is closed under taking direct products, homomorphic images and  subgroups. Moreover, the extension of a
 $\sigma$-soluble group by a $\sigma$-soluble group is  a $\sigma$-soluble group.}

 A  subgroup   $H$  of $G$ is said to be:  a Hall \emph{$\Pi$-subgroup} of $G$ if $|H|$ is $\Pi$-number and $|G:H|$ is $\Pi'$-number;     a \emph{$\sigma$-Hall} subgroup of $G$    if  $H$ is a Hall  $\Pi$-subgroup of $G$ for some $\Pi\subseteq \sigma$.
If $G$ has a complete Hall $\sigma$-set ${\cal H}=\{1, H_{1}, \ldots , H_{t} \}$   such that $H_{i}H_{j}=H_{j}H_{i}$ for all $i, j,$ then we say  that $\{ H_{1}, \ldots , H_{t} \}$ is a \emph{$\sigma$-basis} of $G$.

Let $A$, $B$ and $R$ be subgroups of $G$. Then $A$ is said to \emph{$R$-permute} with  $B$ \cite{guo2} if for some $x\in R$ we have $AB^{x}=B^{x}A$.

The following proposition gives the basic properties of $\sigma$-soluble groups.

{\bf Proposition  3.2.}  {\sl Let $G$ is $\sigma$-soluble. Then:}

(i)  {\sl $|G:M|$ is $\sigma$-primary for every maximal subgroup $M$ of $G$.}

(ii)  {\sl  For every $\sigma _{i}\in \sigma (G)$, $G$ has a  maximal  subgroup $M$ such that $|G:M|$ is a $\sigma _{i}$-number.}

(iii) {\sl  $G$ has  a $\sigma$-basis $\{H_{1}, \ldots , H_{t} \}$ such that for each $i\ne j$ every Sylow subgroup of $H_{i}$  $G$-permutes with every Sylow subgroup of $H_{j}$. }

(iv) {\sl For any $\Pi$, $G$ has a  Hall $\Pi$-subgroup and  every    $\sigma $-Hall subgroup of $G$  $G$-permutes with every Sylow subgroup of $G$.   }

(v) {\sl For any $\Pi$,  $G$ has a  Hall $\Pi$-subgroup $E$,  every  $\Pi$-subgroup of $G$ is contained in some conjugate of $E$ and $E$   $G$-permutes with every
Sylow subgroup of $G$. }

{\bf Proof.}  (i) If $H/M_{G}$ is a chief factor of $G$, then $|(G/M_{G}):(M/M_{G})|=|G:M|$ divides $|H/M_{G}|$, so it is $\sigma$-primary.

(ii) Let $R$ be a minimal normal subgroup of $G$. Then $R$ is a $\sigma _{k}$-group, for some $k$.  If  $R$ is not a Hall $\sigma _{k}$-subgroup of $G$, then  $G/R$ is a  $\sigma$-soluble group such that  $\sigma (G/R)= \sigma (G)$.  Hence by induction, for every  $\sigma _{i}\in \sigma (G/R)$, $G/R$ has a  maximal  subgroup $M/R$ such that
$|(G/R):(M/R)|=|G:M|$ is a $\sigma _{i}$-number. Now suppose that $R$ is  a Hall $\sigma _{k}$-subgroup of $G$ and let $U$ be a complement to $R$ in $G$. Then $G$ has a maximal subgroup $M$ such that $|G:M|$ divides $|R|$, so it is a $\sigma _{k}$-number. On the other hand, for every $\sigma _{i}\ne \sigma _{k}$, $\sigma _{i} \in \sigma(G/R)$ and  so as above we get that $G$ has a maximal subgroup $M$ such that $|G:M|$ is a  $\sigma _{i}$-number.

(iii), (iv), (v) See Theorems A and B in \cite{2}. The proposition is proved.

Let $H/K$ be a chief factor of $G$. Then we say that $H/K$ is \emph{$\sigma$-central} in $G$ if the semidirect product $(H/K) \rtimes (G/C_{G}(H/K))$  is $\sigma$-primary.  Otherwise, it is called \emph{$\sigma$-eccentric } in $G$.

The following lemma is well-known (see for example  Lemma 3.29 in \cite{Shem-Sk}).

{\bf Lemma 3.3.} {\sl Let $R$ be an abelian minimal normal subgroup of $G$ such that $G=RM$ for a  maximal subgroup $M$  of $G$. Then   $G/M_{G}\simeq R\rtimes
 (G/C_{G}(R))$.}

It is well-known that a nilpotent group can be characterized as the group in which each subgroup, or each Sylow subgroup, or each maximal subgroup is subnormal. The following result demonstrates that there is a quite similar relation between  $\sigma$-nilpotency and $\sigma$-subnormality.

{\bf Proposition  3.4.}  {\sl  Any two  of the following conditions are equivalent:  }

(i)  {\sl $G$  is $\sigma $-nilpotent.}

(ii)  {\sl Every chief factor of $G$ is $\sigma$-central in $G$.}

(iii)  {\sl $G$ has a complete Hall $\sigma$-set ${\cal H} $ such that every member of  ${\cal H}$  is $\sigma $-subnormal in $G$.}

(iv)  {\sl Every subgroup of $G$ is ${\sigma}$-subnormal in $G$. }

(v)  {\sl Every maximal subgroup of $G$ is ${\sigma}$-subnormal in $G$.}

{\bf Proof.}  Since (i) $\Rightarrow$ (iii) and (iv) $\Rightarrow$ (v) are clear, it is enough to prove the implications  (i) $\Rightarrow$
(ii), (ii) $\Rightarrow$ (v),   (iii) $\Rightarrow$ (i),  (i) $\Rightarrow$ (iv)  and (v) $\Rightarrow$ (i).

(i) $\Rightarrow$ (ii) For every chief factor $H/K$ of $G$, where $H\leq H_{i}$, we have that $(H/K)\rtimes (G/C_{G}(H/K))$ is a $\pi (H_{i})$-group.  Hence
$(H/K)\rtimes (G/C_{G}(H/K))$ is $\sigma$-primary. Now applaying the Jordan-H\"older theorem \cite[A, 3.2]{DH}, we get that every chief factor of $G$ is $\sigma$-central.

(ii) $\Rightarrow$ (v)  Let $M$ be a maximal subgroup of $G$. Assume that $M_{G}\ne 1$. It is clear that the hypothesis holds for $G/M_{G}$,
 so   $M/M_{G}$ is $\sigma $-subnormal in  $G/M_{G}$ by induction. Hence $M$ is $\sigma $-subnormal in  $G$ by Lemma 2.2(3). Now assume that $M_{G}=1$.
By \cite[A, 15.2]{DH},  either $G$ has a unique minimal normal subgroup $R$ or $G$ has exactly two minimal normal subgroups $R$ and $N$ and the following hold:
  $R$ and $N$ are isomorphic non-abelian groups,  $R\cap M= 1= N\cap M$ and $C_{G}(R)=N$. Let  $C=C_{G}(R)$.  Suppose that $R$ is abelian.    Then $C=R$ by \cite[A,
15.2]{DH}, so in this case we have $G\simeq G/M_{G}\simeq R\rtimes  (G/C_{G}(R))$ is $\sigma$-primary by Lemma 3.3(ii).  Hence, for some
 $\sigma _{i}\in \sigma$, $G$ is a $\sigma_{i}$-group. But then $M$ is $\sigma$-subnormal in $G$. Similarly we get that  $M$ is $\sigma$-subnormal in $G$ in the case when $C=1$.
Finally, if $C=N$,  then $G/N\simeq M\simeq G/R$ is $\sigma$-primary. It follows that $G$ is a $\sigma _{i}$-group. Thus  $M$ is $\sigma$-subnormal in $G$.

(iii) $\Rightarrow$ (i) This follows from    Proposition 2.5 since every member of $\cal H$ is $\sigma$-nilpotent.

(i) $\Rightarrow$ (iv) This follows from the implications (i) $\Rightarrow$ (ii), (ii) $\Rightarrow$ (v) and the evident fact that every subgroup of a $\sigma$-nilpotent group is $\sigma$-nilpotent (see Lemma 2.1).

(v) $\Rightarrow$ (i)    Assume that this  is false and let $G$ be a counterexample of minimal order.

First note that $G$ is $\sigma $-soluble. Indeed, for any maximal subgroup $M$ of $G$, $G/M_{G}$ is $\sigma$-primary and so   $G/M_{G}$ is $\sigma$-soluble. But then $G/\Phi (G)$ is a subdirect product of some $\sigma$-soluble groups, which implies that $G/\Phi (G)$ is $\sigma$-soluble by Lemma 3.1.  Hence   $G$ is  $\sigma$-soluble.  By Proposition 3.2, $G$ has a
a complete Hall $\sigma$-set ${\cal H}=\{1, H_{1}, \ldots , H_{t} \}$.

Let $H=H_{i}$ and $R$ be a minimal normal subgroup of $G$.   We  show that $H$ is normal
in $G$. Assume that is false.   By Lemma 2.2(6), the hypothesis holds for  $G/R$, so $HR/R$ is normal in $G/R$ by the
choice of $G$.  Hence we can assume that $R\nleq H$, so $R\cap H=1$  since $G$ is $\sigma $-soluble. If $G$ has a minimal normal subgroup
$N\ne R$, then as above we get that $HN$ is normal in $G$ and so $RH\cap NH=H(R\cap NH)=H(R\cap N)=H$ is normal in $G$. Therefore $R$ is the unique
minimal normal subgroup of $G$. Moreover, in view of Lemma 2.1, we have  $R\nleq  \Phi (G)$ since $HR/R\simeq H$ is $\sigma$-nilpotent  and  $HR$ is normal in $G$. Let $M$ be a maximal subgroup of $G$ such that
$G=RM$. Then $M_{G}=1$. But $M$ is $\sigma$-subnormal in $G$ by hypothesis and  $G\simeq G/M_{G}$ is $\sigma$-primary, which implies that $G$ is a
$\sigma _{i}$-group, for some $\sigma _{i}\in \sigma$. Therefore $H=G$. This contradiction shows that  (v) $\Rightarrow$ (i).

The proposition is proved.

We say that  $G$ is \emph{$\Pi$-closed} if $O_{\Pi}(G)$ is a Hall $\Pi$-subgroup of $G$.

{\bf Lemma  3.5.}  {\sl Let  $H$ be a normal subgroup of $G$. If
 $H/H\cap \Phi (G)$ is  $\Pi$-closed, then $H$ is $\Pi$-closed}.

{\bf Proof.}  See the proof of   Lemma 2.5 in \cite{1}.

The  integers $n$ and $m$ are called \emph{$\sigma$-coprime} if $\sigma (n) \cap \sigma (m)= \emptyset$.

{\bf Lemma 3.6.}    {\sl  If a $\sigma$-soluble  group $G$ has
   three  $\Pi$-closed  subgroups $A$, $B$ and $C$ whose indices
$|G:A|$, $|G:B|$, $|G:C|$ are pairwise  $\sigma$-coprime, then
$G$ is $\Pi$-closed. }

 {\bf Proof.}  Suppose  that this lemma   is false and let $G$ be
counterexample  of minimal order. Let  $N$ be
 a minimal normal subgroup of $G$. Then  the hypothesis holds
 for $G/N$, so $G/N$ is $\Pi$-closed  by the choice of $G$. Therefore $N$
is not a  $\Pi$-group. Moreover,  $N$ is the unique minimal normal subgroup of $G$
 and, by Lemma 3.5,  $N\nleq \Phi (G)$.
Hence $C_{G}(N)\leq N$.  Since $G$ is   $\sigma$-soluble by
hypothesis, $N$ is  a $\sigma _{i}$-group for some $i$.
Then  $\sigma _{i}\in \Pi'$.

Since $|G:A|$, $|G:B|$, $|G:C|$ are pairwise  $\sigma$-coprime, there are at
least two subgroups, say $A$ and $B$, such  that $N\leq A\cap B$.  Then
$O_{\Pi}(A)\leq C_{G}(N)\leq N$, so  $O_{\Pi}(A)=1$.  But by hypothesis,
$A$ is $\Pi$-closed,  hence  $A$ is a  $\Pi'$-group.
Similarly we get that  $B$ is a  $\Pi'$-group and so $G=AB$ is
 a  $\Pi'$-group. But then $G$ is   $\Pi$-closed. This
contradiction completes the proof of the lemma.

{\bf Lemma 3.7} (See  \cite[III, 5.2]{hupp}
).  {\sl If $G$ is a
 Schmidt group, then $G=P\rtimes Q$, where    $P=G^{\frak{N}}$
 is a Sylow $p$-subgroup of $G$ and $Q=\langle x \rangle $ is a cyclic
 Sylow $q$-subgroup of $G$. Moreover, $\langle x^{q} \rangle \leq \Phi (G)$,   $P/\Phi (P)$ is a chief factor of $G$, $P$
is of exponent $p$ or exponent 4 (if $P$ is a non-abelian 2-group) and $\Phi (P)=1$ if $P$ is abelian.}

{\bf Proposition 3.8.} {\sl  Let $G$ be a $\sigma$-soluble group.
Suppose that $G$    is not $\Pi$-closed but all proper
subgroups of $G$  are $\Pi$-closed.  Then $G$ is  a  $\Pi'$-closed
  Schmidt group.    }

{\bf Proof.}   Suppose that this proposition is false and let $G$ be a
counterexample of minimal order. Let  $R$ be
 a minimal normal subgroup of   $G$ and $\{H_{1}, \ldots ,
H_{t} \}$    a  complete Hall $\sigma$-set of $G$. Without loss of generality we can assume that $H_{i}$ is a $\sigma_{i}$-group for all $i=1, \ldots, t$.

(1)  {\sl   $|\sigma (G)|=2$. Hence $G=H_{1}H_{2}$.}

It is clear that $|\sigma (G)| > 1$. Suppose that $|\sigma (G)| > 2 $.
Then, since $G$ is $\sigma$-soluble, there are maximal subgroups $M_{1}$,
$M_{2}$ and $M_{3}$ whose indices $|G:M_{1}|$, $|G:M_{2}|$  and
$|G:M_{3}|$   are pairwise $\sigma$-coprime. But the  subgroups $M_{1}$,
$M_{2}$ and $M_{3}$ are  $\Pi$-closed   by hypothesis. Hence
 $G$ is $\Pi$-closed  by Lemma 3.6, a contradiction.
Thus $|\sigma (G)|=2$.

Without loss of generality we can assume that  $\sigma_{2}\in \Pi$. Then $\Pi \cap \sigma (G)=\{\sigma _{2}\}$.

(2) {\sl  If either $R\leq \Phi (G)$ or $R\leq H_{2}$, then  $G/R$ is a  $\Pi'$-closed
Schmidt group. }

 Lemma 3.5 and the choice of $G$ imply that  $G/R$ is not $\Pi$-closed.    On the other hand, every
maximal subgroup $M/R$ of $G/R$ is   $\Pi$-closed since $M$ is $\Pi$-closed by hypothesis.
Hence the hypothesis holds for $G/R$. The
choice of $G$   implies that   $G/R$ is a  $\Pi'$-closed
Schmidt group.

(3)  {\sl  $\Phi (G)=1 $,  $R$ is  the unique
minimal normal subgroup of $G$ and $R\leq  H_{1}$.   }

Suppose that    $R\leq \Phi (G)$.  Then  $\bar{G}=G/R$ is a  $\Pi'$-closed
Schmidt group by Claim (2), so  $\bar{G}= \bar{H}_{1}\rtimes \bar{H}_{2}=\bar{P}\rtimes \bar{Q}$,
where  $\bar{H}_{1}=\bar{P}$ is a $p$-group and $\bar{H}_{2}=\bar{Q}$
 is a $q$-group for some primes $p$ and $q$ by Lemma 3.7. Therefore, in fact, $G$  is
not $p$-nilpotent but  every maximal subgroup of $G$ is
$p$-nilpotent.  Hence  $G$ is $\Pi'$-closed
Schmidt group by \cite[IV, 5.4]{hupp}, a contradiction.  Therefore $\Phi (G)=1$.

Now  assume that $G$ has a minimal normal subgroup $L\ne R$.
Since $\Phi (G)=1$, there are maximal subgroups $M$ and $T$   of $G$ such that $LM=G$ and
$RT=G$. By hypothesis, $M$ and $T$ are  $\Pi$-closed. Hence
$G/L\simeq LM/L\simeq M/M\cap L$ is $\Pi$-closed. Similarly,
$G/R$ is $\Pi$-closed and so $G\simeq G/L\cap R$ is $\Pi$-closed, a contradiction.  Hence $R$ is  the unique
minimal normal subgroup of $G$   and  $R\leq  H_{1}$.

  {\sl  Final contradiction.}  In view of Claim (3),
$C_{G}(R) \leq R$ and $R\leq  H_{1}$. Hence $|H_{2}|$ is a prime and
   $RH_{2}=G$ since  every proper subgroup of $G$
 is $\Pi$-closed. Therefore $R=H_{1}$, so $R$ is not abelian since  $G$ is not
a $\Pi'$-closed  Schmidt group.    It is clear that  for any prime
$p$ dividing $|R|$ there is a Sylow $p$-subgroup $P$ of $G$ such that $PH_{2}=H_{2}P$ by Lemma 3.2(iv).
But $H_{2}P < G$, so   $H_{2}P=H_{2}\rtimes  P$. Therefore $R\leq
N_{G}(H_{2})$ and thereby $G=R\times H_{2} =H_{1}\times H_{2}
$ is $\sigma$-nilpotent. This final contradiction completes the proof.

We say that $G$ is \emph{$\sigma$-fiber} if $|\sigma (G)|=|\pi (G)|$.

{\bf Corollary  3.9.}    {\sl Suppose that  $G$ is not $\sigma$-nilpotent but every proper subgroup of $G$ is $\sigma$-nilpotent. If $G$ is  $\sigma$-soluble, then $G$ is a $\sigma$-fiber Schmidt group.}

{\bf Proof.} It is clear that $G$ is $\sigma$-nilpotent if and only if $G$ is $\Pi$-closed for all $\Pi\subseteq \sigma$. Hence, for some $\Pi$, $G$
is not $\Pi$-closed. On the other hand, every proper subgroup of $G$ is $\Pi$-closed. Hence  $G$ is a Schmidt group by Proposition 3.8 and clearly $|\sigma (G)|=\pi(G)|.$

\section{Proof of Theorem 1.2}

In this section, we  need the following

{\bf Lemma 4.1}  (See Lemmas 2.8, 3.1 and Theorem B(1) in \cite{1}).  {\sl  Let  $H$,  $K$ and $R$ be subgroups of a $\sigma$-soluble group $G$. Suppose that  $H$ is $\sigma$-quasinormal in $G$ and $R$ is normal in $G$. Then:}

(1)  {\sl If  $H\leq E \leq G$, then $H$ is $\sigma$-quasinormal in $E$.}

(2)  {\sl The subgroup  $HR/R$ is $\sigma$-quasinormal in $G/R$.}

(3) {\sl If  $R\leq H$ and $H/R$ is $\sigma$-quasinormal in $G/R$, then $H$ is  $\sigma$-quasinormal in $G$.}

(4) {\sl $H$ is $\sigma$-subnormal in $G$. }

(5) {\sl If $H$ is a $\sigma _{i}$-group, then $O^{\sigma _{i}}(G)  \leq N_{G}(H)$ }.

{\bf Lemma  4.2.}  {\sl The following statements hold:}

(1) {\sl  $m_{\sigma}(G)\leq m_{\sigma q}(G)$.}

(2) {\sl  If $M$ is a non-$\sigma$-subnormal maximal subgroup of $G$, then $m_{\sigma}(M)\leq m_{\sigma}(G)-1$.}

(3) {\sl If $R$ is a normal subgroup of $G$, then  $m_{\sigma}(G/R)\leq m_{\sigma}(G)$.}

(4) {\sl If $R$ is a normal subgroup of $G$ and $G$ is $\sigma$-soluble, then  $m_{\sigma q}(G/R)\leq m_{\sigma q}(G)$.}

{\bf Proof.} (1) This follows from   Lemma   4.1(4).

(2) Since $M$ is not $\sigma$-subnormal in  $G$,  $m_{\sigma q}(G) > 1.$
Moreover, for $n=m_{\sigma }(G)$, each $(n-1)$-maximal subgroup of $M$ is $\sigma$-subnormal in $M$ by Lemma 2.2(1). Hence $m_{\sigma}(M)\leq m_{\sigma}(G)-1$.

(3) If  each  maximal chain of $G/R$ has  length $r < m_{\sigma}(G)$,   it is  clear. Otherwise, this follows from Lemma 2.2(3)(6).

(4) This is a corollary of  Lemma 4.1(3)

The lemma is proved.

The following  properties of the rank of a soluble group are useful in our proof.

{\bf  Lemma 4.3}   (See \cite[VI, Lemma 5.3]{hupp}). {\sl Let $G$ be soluble. Then:}

(1) {\sl $r(G/R)\leq r(G)$ for all  normal subgroups $R$ of $G$}

(2) {\sl $r(E)\leq r(G)$ for all  subgroups $E$ of $G$ }

(3) $r(A\times B)= \text{Max} \{r(A), r(B)\}$.

{\bf Lemma   4.4} (See Huppert \cite[Lemma 11]{HupI}). {\sl Let $G$ be soluble and $R$ be a minimal normal subgroup of $G$. Let $H$ be a minimal supplement to $C_{G}(R)$ in $G$. Then $H\cap R=1$.}

{\bf Lemma 4.5}. {\sl The following statements hold:}

 (i) {\sl If each $n$-maximal subgroup of $G$ is $\sigma$-subnormal and $n > 1$, then each $(n-1)$-maximal subgroup is $\sigma$-nilpotent.}

(ii) {\sl  If each $n$-maximal subgroup of $G$ is $\sigma$-subnormal, then each $(n+1)$-maximal subgroup is $\sigma$-subnormal.}

{\bf Proof.}    (i) Let $H$ be an $(n-1)$-maximal subgroup of $G$ and $K$  a maximal subgroup of $H$. Then $K$ is  an $n$-maximal subgroup of $G$, so it is $\sigma$-subnormal in $G$. Then, by Lemma 2.2(1), $K$ is $\sigma$-subnormal in $H$. Therefore each maximal subgroup of $H$ is
$\sigma$-subnormal in $H$. It follows from Proposition 3.4 that $H$ is $\sigma$-nilpotent.

(ii) Let $L \leq M\leq G$, where  $M$ is an $n$-maximal subgroup of $G$ and $L$ is a maximal subgroup of $M$. If $n=1$, $G$ is $\sigma$-nilpotent and  so $L$ is $\sigma$-subnormal in $G$  by Proposition 3.4. On the other hand, in the case when  $n > 1$ Statement (i) implies that
  each $(n-1)$-maximal subgroup of $G$ in $\sigma$-nilpotent. Then $M$ is $\sigma$-nilpotent by Lemma 2.1, so $L$ is $\sigma$-subnormal in $G$ by
  Lemma 2.2(4).

The lemma is proved.

{\bf  Proof of Theorem 1.2.}  Suppose that this theorem  is false and let $G$ be a counterexample of minimal order. Let   ${\cal H}  =\{1, H_{1}, \ldots , H_{t}\}$ be a complete Hall $\sigma$-set of $G$. Then $t > 1$.

(i)  Suppose that this   is false. Let $R$ be a minimal normal subgroup of $G$ and $|R|=p^{m}$. Without loss of generality we can assume that $R\leq H_{1}$.
 Let   $n=m_{\sigma q}(G)$. Since  $G$ is not $\sigma$-nilpotent, some maximal subgroup $M$  of $G$ is not $\sigma$-subnormal in $G$ by Proposition 3.4 and  so $M$ is not $\sigma$-quasinormal in $G$ by Lemma 4.1(4). Thus $n > 1$.

(1)  {\sl   $r(G/R) \leq n+r-2$. }

Assume that  $r(G/R) > n+ r -2$.  Note that $\{H_{1}R/R, \ldots , H_{t}R/R\}$ is  a complete Hall
 $\sigma$-set of $G/R$  and   $r(H_{i}R/R) = r(H_{i}/H_{i}\cap R)\leq r(H_{i})\leq r$ for all $i=1, \ldots, t$ by Lemma 4.3(1). Assume that $G/R$ is $\sigma$-nilpotent.
Then $G/R=  (H_{1}R/R) \times \cdots \times  (H_{t}R/R)$, so  $r(G/R)  \leq r\leq  n+r-2$  since $n > 1$ by Lemma 4.3(3). This contradiction shows that
$G/R$ is not $\sigma$-nilpotent. Moreover, $G/R$ is $\sigma$-soluble  by Lemma 3.1. By Lemma 4.1(2)(3), $m_{\sigma q}(G/R)\leq m_{\sigma q}(G)=n.$ The Choice of $G$ implies that
$r(G/R) \leq m_{\sigma q}(G/R)+r-2\leq  n+r -2$, a
 contradiction.    Hence we have (1).

(2) {\sl   $m >  n+ r -2$. Hence $R$  is the only minimal normal subgroup of $G$. }

First note that in view of the Jordan-H\"{o}lder theorem, Claim (1) and the choice of
 $G$  we have $m > n + r-2$.
If $G$ has a minimal normal subgroup $N\ne R$, then
 $r(G/N) \leq n + r-2$ by Claim (1), so in view of the $G$-isomorphism
 $R\simeq RN/N$ we get that
 $m \leq n + r-2$, a contradiction. Hence   $R$
 is the only minimal normal subgroup of $G$.

(3) {\sl If $M$ is a proper subgroup of $G$, then $r(M) \leq n +  r-2$.}

It is enough to consider the case when $M$ is a maximal subgroup of $G$.
Assume that  $r(M) > n +   r-2$. Then $M$ is not $\sigma$-nilpotent (see the proof of  Claim (1)).
 Therefore $n > 2$ by Lemmas 4.1(1)  and Proposition 3.4.
Moreover, since  $G$ is $\sigma$-soluble,  $M$ possesses a complete Hall
$\sigma$-set $\{M_{1}, \ldots, M_{t}\}$ such that $M_{i}=H_{i}\cap M$
for all $i=1, \ldots , t$ by Lemma 3.2(v).  Hence  $r(M_{i})\leq r(H_{i})  \leq r$ for all   $i=1, \ldots , t$ by Lemma 4.3(2).
Therefore,  $M$ satisfies the hypothesis, with $n-1$ instead of $n$,
   by Lemmas 4.1(1) and  so
 the choice of $G$ implies that
 $r(M) \leq  n-1 + r -2  \leq n +  r-2$, a contradiction.  Hence we have (3).

(4) {\sl $R\nleq \Phi (G)$.}

Suppose that  $R\leq \Phi (G)$. Then for a minimal supplement $H$ to $C_{G}(R)$
 in $G$ we have $H\cap R=1$ by Lemma 4.4, so $RH\ne G$ and $R$ is a minimal
 normal subgroup of $RH$. But Claim (3) implies that $r(RH) \leq n + r-2$ and
  so $m\leq n + r-2$, contrary to Claim (2). Hence we have (4).

{\sl Final contradiction for (i). }  Claim (4) implies  that there is a
 maximal subgroup $M$ of $G$ such that
  $G=RM$ and $H_{2}\leq M$. Then  $M_{G}=1$ by Claim (2), so
 $C_{G}(R)=C_{G}(R)\cap RM=R(C_{G}(R) \cap M)=R$.
Let $H_{2}= M_{s}$  be a member of a maximal chain $1=M_{l} < M_{l-1}
 < \ldots < M_{1} < M_{0}=M$ of $M$.
Then $R$ is an  $l$-maximal subgroup of $G$.  First suppose that $l > n-1$.
Assume also that $H_{2} \leq M_{n-1}$.  By hypothesis $M_{n-1}$ is
$\sigma$-quasinormal in $G$. Hence $H^{x}_{2}\leq M_{n-1}$ for all $x\in
G$. It follows that $(H_{2})^{G}\leq M_{G}=1$, so
$H_{2}=1$, a contradiction. Therefore  $n\leq s$,
 that is, for the $n$-maximal subgroup $H=M_{n-1}$ of $G$
contained in $H_{2}$   we have $H\ne 1$.
  Then   $H\leq O_{\pi (H_{2})}(G)$ by Lemma 2.2(9).
But   $R\cap O_{\pi (H_{2})}(G)=1$ since $R\leq H_{1}$,  so  $1 < H \leq
O_{\pi (H_{2})}(G)\leq C_{G}(R)=R$,  a   contradiction.

Therefore $n-1 \leq l$, so $M$  possesses a  maximal chain  $1=M_{k} < M_{k-1}
 < \cdots  M_{1} <  M_{0}=M$, where $k < n$.
Then   $R$ is a $k$-maximal subgroup of $G$. Therefore every $l$-maximal subgroup of $R$ is a
$(k+l)$-maximal subgroup of $G$.
 Let $R_{0}$ be a minimal normal
 subgroup of $H_{1}$ contained in $R$ with $|R_{0}|=p^{a}$.
  Let   $L$ be   an $(n-k)$-maximal subgroup
 of $R$  with    $|L|=p^{b}$  such that $L\leq R_{0}$ in the case when $b < a$
and   $R_{0}\leq L$ if $a\leq b$.  Then $L$ is    an $n$-maximal subgroup
 of $G$, so $L$ is $\sigma$-quasinormal in $G$.

First suppose that  $L\leq R_{0}$. Then $$L^{G}
 =L^{H_{1}N_{G}(L)}=L^{H_{1}}\leq (R_{0})_{G}=R$$ by Lemma 4.1(5). Hence
$R_{0}=R$.  Then $m=a\leq r$ and  so $m\leq r+n-2$ since $n > 1$, contrary to (2).
 Thus   $R_0\leq L$, so  $$R_{0}^{G} =R_{0}^{H_{1}N_{G}(L)}=R_{0}^{N_{G}(L)}\leq L,$$
which implies that $L=R$, a contradiction also.  Hence Assertion (i) is true.

 (ii) Let $n=m_{\sigma}(G)$.  Suppose that $n > l_{\sigma}(G) $. Then $n  > 1$. Indeed, if $n=1$, $G$ is $\sigma$-nilpotent by Proposition 3.4
  and so $l_{\sigma}(G)=1 =m_{\sigma}(G)$, a contradiction.
  The choice of $G$ and Lemma 4.2(3) imply that $l_{\sigma}(G/F_{\sigma}(G))\leq m_{\sigma}(G/F_{\sigma}(G))\leq  n$.
  Hence  $F_{\sigma}(G)\nleq \Phi (G)$ by Lemma 2.3(2). Let $M$ be a maximal subgroup of $G$ such that $F_{\sigma}(G)M=G$.
   Then $l_{\sigma}(G/F_{\sigma}(G))=l_{\sigma}(MF_{\sigma}(G)/F_{\sigma}(G)=l_{\sigma}(M/M\cap F_{\sigma}(G))=l_{\sigma}(G)-1$.  Lemma 2.1 implies that $M\cap F_{\sigma}(G)\leq F_{\sigma}(M)$. Therefore  $l_{\sigma}(G) \leq l_{\sigma}(M)+1$.
   Note that $m_{\sigma}(M)  \leq   m_{\sigma}(G) -
   1$. Indeed, since  $ n > 1$,  each $(n-1)$-maximal subgroup of $M$ is $\sigma$-quasinormal in $G$. Hence each $(n-1)$-maximal subgroup of $M$ is $\sigma$-quasinormal in $M$ by Lemma 4.1(1). Hence $m_{\sigma}(M)\leq n-1$. But the choice of $G$ we have $ l_{\sigma}(M)\leq  m_{\sigma}(M)$,  and then $ l_{\sigma}(G)\leq  l_{\sigma}(M) +1 \leq  m_{\sigma}(M) +1 \leq m_{\sigma}(G)$. This contradiction completes the proof of (ii).

 (iii)  Suppose that   $ m_{\sigma}(G)  < |\pi (G)|$. Let $P_{1}, \ldots , P_{n}$ be a Sylow basis of $G$  and  $\cal H$ a complete Hall $\sigma$-set of $G$. Then for any $i$, $i=1$ say, we have
 $P_{1} < P_{1}P_{2}< \cdots < P_{1}P_{2}\cdots P_{n}=G$, so $P_{1}$ is at least an $(n-1)$-maximal subgroup of $G$. Therefore $P_{1}$ is $\sigma$-subnormal in $G$ by Lemma 4.5(ii) since $ m_{\sigma}(G)  < |\pi (G)|$. Hence every Sylow subgroup of $G$ is $\sigma$-subnormal in $G$ and so every member of $\cal H$ is $\sigma$-subnormal in $G$ by Proposition 2.5. But then $G$ is $\sigma$-nilpotent by Proposition 3.4, a contradiction. Hence $|\pi (G)|\leq m_{\sigma}(G)$.

 The theorem  is proved.

\section{Proofs of Theorem 1.4  and Corollaries 1.7 and 1.8}

{\bf Proof of Theorem 1.4. } Let $R$ be a minimal normal subgroup of $G$.

(i) Suppose that this assertion  is false and let $G$ be a counterexample of minimal order.
  First note   that  $G/R$  is
$\sigma$-soluble. Indeed, if $R$ is a maximal subgroup or a 2-maximal
subgroup of $G$, it is clear. Otherwise,   the hypothesis
holds for $G/R$ by  Lemma 2.2(6), so the choice of $G$ implies that  $G/R$ is
$\sigma$-soluble. Hence $R$ is the unique minimal normal subgroup of $G$ by Lemma 3.1  and $R$ is not $\sigma$-primary. Hence $R$ is not abelian.

 Let $p$ be any odd prime dividing $|R|$ and $R_{p}$ a Sylow $p$-subgroup of $R$.
 The Frattini argument  implies  that there is a maximal subgroup  $M$
 of $G$ such that $N_{G}(R_{p})\leq M$ and $G=RM$. It is clear that $M_{G}=1$, so $M$
  is not $\sigma$-subnormal in $G$ since $G/M_{G}\simeq G$ is not $\sigma$-primary.
   Let $D=M\cap R$. Then $R_{p}$ is a Sylow $p$-subgroup of $D$.

(1) {\sl $D$ is not nilpotent.   Hence $D\nleq  \Phi (M)$ and $D$ is not a $p$-group.}

Assume that  $D$ is a nilpotent. Then
 $R_{p}$ is normal in $M$. Hence $Z(J(R_{p}))$ is normal
 in $M$. Since $M_{G}=1$, it follows that $N_{G}(Z(J(R_{p})))=M$ and so
 $N_{R}(Z(J(R_{p})))=D$   is  nilpotent. This implies that $R$ is $p$-nilpotent
 by Glauberman-Thompson's theorem on
 the normal $p$-complements. But
then $R$ is a $p$-group, a contradiction. Hence we have (1).

(2) $R < G$.

Suppose that $R=G$ is a simple non-abelian group. Assume that some proper
 non-identity subgroup $A$ of $G$ is $\sigma$-subnormal in $G$. Then there
 is a subgroup chain
 $A=A_{0} \leq A_{1} \leq \cdots \leq A_{n}=G$  such that  either $A_{i-1}\trianglelefteq  A_{i}$ or
$A_{i}/(A_{i-1})_{A_{i}}$ is  ${\sigma}$-primary for all $i=1, \ldots ,t$.
  Without loss of generality, we can assume that $M=A_{n-1} < G$. Then $M_{G}=1$
 since  $G=R$ is simple, so $G\simeq G/1$ is $\sigma$-primary,
 a contradiction. Hence every proper $\sigma$-subnormal
 subgroup of $G$ is trivial.

Let $P$ be a Sylow $p$-subgroup of $G$,  where $p$ is the smallest prime dividing $|G|$,
  and let $L$ be a maximal subgroup  of $G$ containing $P$.
Then, in view of \cite[IV, 2.8]{hupp}, $|P| > p$.  Let $V$ be a maximal subgroup of $P$.
 If $|V|=p$, then $P$ is abelian, so $P < L$ by \cite[IV, 7.4]{hupp}. Hence there is a 3-maximal subgroup $W$ of
 $G$ such that $V\leq W$. But then some proper non-identity  subgroup of $G$ is
 $\sigma$-subnormal in $G$ by hypothesis, a contradiction. Therefore $|V| > p$, which again
 implies that some proper non-identity  subgroup of $G$ is $\sigma$-subnormal
 in $G$.  This contradiction shows that we have (2).

(3) {\sl $M$ is $\sigma$-soluble.}

Let $L < T < M$, where $L$ is a maximal subgroup of $T$ and $T$ is a maximal subgroup of $M$.
Since $M$ is not $\sigma$-subnormal in $G$, either $L$ or $T$ is
 $\sigma$-subnormal in $G$ and so it is $\sigma$-subnormal in $M$  by
Lemma 2.2(1). Hence the hypothesis holds for $M$,  so $M$ is $\sigma$-soluble
 by the choice of $G$.

(4) {\sl $M=D\rtimes T$, where $T$ is a maximal   subgroup of $M$ of  prime order.}

In view of Claim (1), there is a maximal subgroup $T$ of $M$ such that $M=DT$.
 Then $G=RM=R(DT)=RT$ and so,
in view of (2), $T\ne 1$.  Assume that $|T|$ is not a prime and let $V$ be a
 maximal subgroup of $T$. Since $M$ is not $\sigma$-subnormal in $G$, at least one of the subgroups $T$ or $V$ is  $\sigma$-subnormal in $G$ by hypothesis.
Claim (3) implies that     $V$ and $T$ are  $\sigma$-soluble.
 Consider, for example, the case when $V$
 is $\sigma$-subnormal in $G$. Since $V\ne 1$ and $V$ is $\sigma$-soluble,
for some $i$ we have $O_{\sigma _{i}}(V)\ne 1$. But
 $O_{\sigma _{i}}(V)\leq O_{\sigma _{i}}(G)$   by Lemma 2.2(9), so $R$ is
$\sigma$-primary, a contradiction.
 Hence $|T|$ is a prime, so $M=D\rtimes T$.

{\sl Final  contradiction for (i). } Since $T$ is a maximal subgroup of $M$ and
 it is cyclic, $M$ is soluble and so $|D|$ is a prime power, which contradicts (1).
 Hence Assertion (i) is true.

(ii)  Suppose that this false. Then   $2\in \pi (G)$.  Part (i) implies that $G$ is $\sigma$-soluble.
Let ${\cal H}=\{H_{1},   \ldots, H_{t}\}$ be a $\sigma$-basis of $G$. Without loss of generality we can assume that $H_{1}$ is a $\sigma _{1}$-group and $2\in \pi (H_{1})$. Then $H_{1}$ is   not soluble,  so $|\pi (H_{1}) | > 1$. Let $p\in \pi (H_{1})$.

(1) {\sl $t=2$ and $H_{2}$ is a Sylow subgroup of $G$}.

By Proposition 3.2(v), $G$ has a  Hall $\sigma _{1} '$-subgroup $E$ and $E$ permutes with some Sylow $p$-subgroup $P$ of $G$ for each $p\in \pi (H_1).$ It is clear that $EP<G$. We show that $PE$ is soluble. In fact, if $PE$ is $\sigma$-nilpotent, then $PE=P\times E$, where $2 \nmid |E|$, so $PE$ is soluble. Now assume that $PE$ is not $\sigma$-nilpotent. Then the hypothesis holds for $PE$, so $PE$ is soluble by the choice of $G$. Hence $PE$ has a Sylow basis  ${\cal P}=\{P, P_{1}, \ldots P_{n}\}$.
If $t>2$ or $H_2$ is not a Sylow subgroup of $G$, then every member of ${\cal P}$ is at least $3$-maximal subgroup of $G$. Hence
every member of ${\cal P}$  is   $\sigma $-subnormal in $G$ by Lemma  4.5(ii). This shows that every Sylow subgroup of $G$ is
 $\sigma $-subnormal in $G$. Therefore all members of $\cal H$ are normal in $G$ by   Lemma 2.2(9), which implies that $G$ is $\sigma$-nilpotent. This contradiction shows that we have (1).

(2)  $O_{\sigma _{1}}(G)\ne 1$.

Suppose that this is false. Since $G$ is $\sigma$-soluble, $R$ is $\sigma$-primary. Hence $R\leq H_{2}$. Let $P$ be a Sylow $p$-subgroup of $H_{1}$ and $P\leq M$, where $M$ is a maximal subgroup of $H_{1}$. Since $H_{1}$ is not soluble, $P < M < H_{1} < G$ by \cite[IV, 7.4]{hupp}.
Therefore there is a 3-maximal subgroup $W$ of $G$ such that $P \leq W\leq H_{1}$. Then $W$ is $\sigma$-subnormal in $G$, so $1 < W\leq O_{\sigma _{1}}(G)$ by Lemma 2.2(9), a contradiction.

Hence we can assume that $R\leq O_{\sigma _{1}}(G)$.

(3) {\sl $H_{2}$ is normal in $G$.}

First suppose that $R$ is a $p$-group for some prime $p$ and let $Q$ be a Sylow $q$-subgroup of $H_{1}$, where $q\ne p$. By Proposition 3.2,  there is $x\in G$ such that $H_{2}Q^{x}=Q^{x}H_{2}$. Hence we have a subgroup chain $H_{2} < H_{2}Q^{x} < RH_{2}Q^{x} < G$. It follows from Lemma  4.5(ii) that $H_{2}$ is $\sigma$-subnormal in $G$, so it is normal in $G$ by Lemma 2.2(9).

Now assume that $R$ is not abelian. Then, for any odd prime $p$ dividing $|R|$, the subgroup $R$   is not $p$-nilpotent. Hence by
Glauberman-Thompson's theorem on the normal $p$-complement, we have that $P< N_{R}(Z(J(P))) < R$, where $P$ is a Sylow $p$-subgroup of $R$.
Since $(|R|, H_{2}))=1$, $H_{2}$ normalizes some Sylow $p$-subgroup of $R$, say $P$. Hence $H_{2}\leq  N_{G}(Z(J(P)))$.
But then we have a chain $H_{2} <   H_{2}P  <   H_{2}N_{R}(Z(J(P))))  < G$, so  $H_{2}$ is $\sigma$-subnormal in $G$ by Lemma 4.5(ii). Consequently  $H_{2}$  is normal in $G$ by Lemma 2.2(9).

{\sl Final contradiction for (ii) }. Let $P$ be a Sylow $p$-subgroup of $H_{1}$ and $V$ be a maximal subgroup of $P$. Since $|\pi (|H_{1}|)| > 1$, $V$ is $\sigma$-subnormal in $G$ by Lemma 4.5(ii). Suppose that $P\nleq O_{\sigma _{1}}(G)$. Then $P\nleq R$  and   $P$ is not $\sigma$-subnormal in $G$ by Lemma 2.2(9), so  $P$ is cyclic by Proposition 2.5. Hence $R\cap P\leq \Phi (R)$, and so $R$ is $p$-nilpotent by the Tate theorem \cite[IV, 4.7]{hupp}. But then $R$ is a $p$-group or a $p'$-group.
Assume that $R$ is a $p$-group. Then $G/R$ is not $\sigma$-nilpotent. Otherwise, $H_1/R\lhd G/R.$ It follows from Claim (3) that $G=H_1\times H_2$, which contradicts $m_{\sigma} (G)>1.$ Hence $G/R$ is not $\sigma$-nilpotent, and so $1< m_{\sigma}(G/R).$ But then $G/R$ satisfies the hypothesis by Lemma 2.2(6). Hence $G/R$ is soluble by induction. Consequently, $G$ is soluble, a contradiction.
Now assume that $R$ is a $p'$-group. Then by the  Frittini argument, $P$ normalizes some Sylow $q$-subgroup $Q$ of $R$, where $q\ne p$ divides $R$. Hence we have a subgroup chain $P < PQ < H_{2}PQ < G$, so $P$ is $\sigma$-subnormal in $G$ by Lemma 4.5(ii). This shows that every Sylow subgroup of $H_{1}$ is  $\sigma$-subnormal in $G$. It follows Lemma 2.2(9) that $H_{1}$ is  $\sigma$-subnormal in $G$. Thus $H_{1}$ is normal in $G$ by Lemma 2.2(9), so $G=H_{1}\times H_{2}$ is $\sigma$-nilpotent. It follows from Proposition 3.4 that $m_{\sigma }(G)=1.$   This final contradiction completes the proof of (ii).

The theorem is proved.

{\bf Proof of Corollary 1.7.} The sufficiency is clear. We only need to prove the necessity.  By Theorem 1.4, $G$ is $\sigma$-soluble. Since $m_{\sigma q}(G)=2$, $G$ is not $\sigma$-nilpotent by Proposition 3.4.
On the other hand, if $M$ is a  maximal subgroup of $G$, then every maximal subgroup of $M$ is $\sigma$-subnormal in $G$ and so it is
$\sigma$-subnormal in $M$ by Lemma 2.2(1). Therefore $M$ is $\sigma$-nilpotent by Proposition 3.4.
Hence $G$ is a Schmidt group such that $|\pi (G)|=|\sigma (G)|$ by Corollary 3.9.
Then by Lemma 3.7,
 $G=P\rtimes Q$, where    $P=G^{\frak{N}}$
 is a Sylow $p$-subgroup of $G$ and $Q=\langle x \rangle $ is a cyclic
 Sylow $q$-subgroup of $G$. Moreover,  $P$
is of exponent $p$ or exponent 4 (if $P$
 is a non-abelian 2-group)   and $P/\Phi (P)$ is a chief factor of $G$.
 If $\Phi (P)\neq 1$, then there exists a maximal subgroup $M$ of $G$ such that $Q< M <G$. By the hypothesis, Then  $Q$ is $\sigma$-subnormal in $G$ by Lemma 4.5(ii). It follows from Lemma 2.2(9) that $Q$ is normal in $G$. This contradiction shows that every Sylow subgroup of $G$ is abelian.

{\bf Proof of Corollary 1.8.} In view of Corollary 1.7 and Lemma 3.7, $G=P\rtimes Q$ is a Schmidt group with $|\pi (G)|=|\sigma (G)|$, where $P$ is a minimal normal subgroup of $G$ and $Q$ is cyclic. Let  $|P| =p^{n}$ and $|Q|=q^{m}$. Suppose that $n > 1$.  Then $G$ has a 2-maximal subgroup $L$ such that $|G:L|=pq$. By hypothesis $L$ is $\sigma$-quasinormal in $G$, so it is, in fact, $S$-quasinormal in $G$ and hence $LQ=QL$ is a subgroup of $G$ with $G:LQ|=p$. But then $LQ\cap P$ is normal in $LQ$  and  $|P:(LQ\cap P)|=p$, so $LQ\cap P$ is normal in $G$ and hence $LQ\cap P=1$ in view of the minimality of $P$. It follows that $|P|=p$, a contradiction. Hence $|P|=p$, so $G$ is supersoluble. The corollary is proved.

\section{Proof of Theorem 1.10}

{\bf Lemma  6.1.} {\sl  Suppose that $G$ is $\sigma $-soluble and let  ${\cal H}= \{H_{1}, \ldots, H_{t}\}$ be a  $\sigma $-basis of $G$. If $H_{i}$ forms an   irreducible pair  with $H_{j}$, then $H_{j}$ is an elementary abelian  Sylow subgroup of $G$.  }

{\bf Proof.} Without loss of generality we can assume that $G=H_{i}H_{j}$. By Proposition 3.2(iv), for each prime $p$ dividing $|H_{j}|$,   there is a Sylow $p$-subgroup $P$ of $H_{j}$ such that $H_{i}P=PH_{i}$. Hence $G=H_{i}P$. Let $R$ be a minimal normal subgroup of $G$. If $R\leq P$, then the maximality of $H_{i}$ implies that $R=P$ is elementary. On the other hand, if $R\leq H_{i}$, then $H_{i}/R $ is a maximal subgroup of $G/R$, and so $P\simeq PR/R$ is elementary by induction.

The following lemma can be proved similarly to \cite[I, Proposition 4.16]{DH}.

{\bf Lemma   6.2.} {\sl Suppose that $G\ne 1$ is $\sigma$-soluble and let $L$ be a subgroup of $G$. Then for each $\sigma$-basis ${\cal L}= \{L_{1}, \ldots, L_{r}\}$   of $L$, there is a $\sigma$-basis  ${\cal H}= \{H_{1}, \ldots, H_{t}\}$   of $G$ such that $L_{i}=L\cap H_{i}$ for all $i=1, \ldots r$.}

{\bf Lemma  6.3.} {\sl Suppose that $G$ is $\sigma$-soluble and let  ${\cal H}= \{H_{1}, \ldots, H_{n}\}$ be a  $\sigma$-basis of $G$. If $H_{1}$ forms an irreducible pair with    $H_{i}^{x}$ for all $i > 1$ and  $x\in G$ such that $H_1H^x_i=H^x_iH_1$,  then  every  subgroup $K$ of $G$ containing $H_{1}$ is a $\sigma$-Hall subgroup of $G$.}

{\bf Proof.} Suppose that this is false.    Without   loss of generality we can assume that $H_{i}$  is a $\sigma _{i}$-group. Let ${\cal K}= \{
H_{1}, K_{2}, \ldots, K_{r}\}$ be a  $\sigma$-basis of $K$. By Lemma 6.2
, there is a $\sigma$-basis $ \{H_{1}, H_{2}^{x_{2}}, \ldots, H_{n}^{x_{n}}\}$ of $G$ such that $K_{i}=K\cap H_{i}^{x_{i}}$ for all $i=2, \ldots , r$.   Hence   $H_{1}  H_{j}^{x_{j}}\cap  K=H_{1}(H_{j}^{x_{j}}\cap K)=H_{1}K_{j}=K_{j}H_{1}$ and so  there is  a subgroup chain  $H_{1} \leq H_{1}K_{j} \leq  H_{1}  H_{j}^{x_{j}}$. The maximality of $H_{1}$ in $H_{1}  H_{j}^{x_{j}}$ implies that $K_{j}=H_{j}^{x_{j}}$. Thus  $K$  is a $\sigma$-Hall subgroup of $G$.

{\bf Lemma   6.4.} {\sl Suppose that $G$ is $\sigma$-soluble and let $K$   be a subgroup of $G$. If every subgroup of  $G$  containing  $K$ is  a
$\sigma$-Hall subgroup of $G$, then $K$ is a $k$-maximal subgroup of $G$, where  $k=|\sigma (|G:K|)|$, and $K$ is not a
 $r$-maximal subgroup of $G$ for all $r > k$. }

{\bf Proof.} Let $\{H_{1}, \ldots , H_{t}\}$ be a $\sigma$-basis of $G$.  The assertion follows from the  fact that  in any maximal chain $K=M_k < M_{k-1} < \ldots < M_1 < M_{0}=G$, $ |M_{i}:M_{i+1}|$ is an order of some $H_{i}$ since  both  $M_{i}$ and $M_{i+1}$ are $\sigma$-Hall subgroups of $G$. The lemma is proved.

{\bf Proof of Theorem 1.10.}   Let $P_{1},  \ldots , P_{n}$ be  a  Sylow  basis of $G$ and
 $\{H_{1}, \ldots , H_{t}\}$  a complete Hall $\sigma$-set of $G$
 We can assume without loss of generality  that each $P_{i}$ is contained in some $H_{j}$ and $P_{i}$ is a $p_{i}$-group.

 {\sl Necessity.} First note that if $G$ is $\sigma$-nilpotent, then $m_{\sigma}(G)=|\pi (G)|=1$, and so $G$ is a $p$-group for some prime $p$. Now we show, assuming that $G$ is not $\sigma$-nilpotent, then $G$ is a  group of type (ii).

 Assume  that this is false and let $G$ be a counterexample of minimal order.  Let $R$ be a minimal normal subgroup of $G$.
   Without loss of generality we can assume that $ R\leq H_{1}$ and  $H_{k}$ is a $\sigma _{i_{k}}$-group for all $k=1, \ldots , t$.

(1) {\sl If $G/R$ is not $\sigma$-nilpotent, then  $G/R$  is  a group of the type (ii).}

Suppose that $G/R$ is not $\sigma$-nilpotent. We show that  the hypothesis holds for $G/R$. Indeed, if $R <   P_{i} $, it is clear. We may, therefore, assume  that  $R =   P_{i} $. Then $R$ has a complement $M$ in $G$ such that $G=R\rtimes M$. Since  $|\pi (M)|=n-1$, $M$ satisfies the same
assumptions as $G$, with $n-1$ replacing $n$, by Lemma 2.2(1). The choice of $G$ implies that  $G/R\simeq M$ is a group of the type (ii).

(2) {\sl If $V_{i}$ is a maximal subgroup of $P_{i}$, then $V_{i}$ is $\sigma$-subnormal in $G$. Hence every non-$\sigma$-subnormal Sylow subgroup of  $G$ is cyclic.}

Since $P_{1},  \ldots , P_{n}$ is  a  Sylow  basis of $G$, $V_{i}$ is at $m$-maximal subgroup of $G$, where $m>n$. Hence
$V_{i}$ is $\sigma$-subnormal in $G$ by Lemma 4.5(ii) Therefore, if $P_{i}$ is not $\sigma$-subnormal in $G$, then it is   cyclic  by Proposition 2.5.

(3) {\sl If $R$ is  the only minimal normal subgroup of $G$, then each Sylow subgroup $P_{i}$ of $H_{k}$ has prime
 order  and it is not $\sigma$-subnormal in $G$   for all $k > 1$.}

Indeed, let $V$ be a maximal subgroup of $P_{i}$. Then $V$ is $\sigma$-subnormal in $G$ by Claim (2). Hence $V\leq O_{\sigma _{i_{k}}}(G)$  by Lemma 2.2(9).
But since $R\leq H_{1}$, we have that $O_{\sigma _{i_{k}}}(G)=1$ and so $V=1$. Moreover, if $P_i$ is $\sigma$-subnormal, then $P_i\leq O_{\sigma _{i_k}}(G)=1,$ a contradiction. Hence we have (3).

(4)  {\sl For some $i$,  $i=1$ say, $P_{i}=P_{1}$ is not $\sigma$-subnormal in $G$. Hence  $P_{1}$  forms an irreducible pair with $P_{i}$ for all $i > 1$.  }

If  $P_{i}$ is    $\sigma$-subnormal in $G$ for all  $i=1, \ldots , n$, then $H_{k}$ is normal in $G$ for all $k=1, \ldots, t$ by Lemma 2.2(9)  and Proposition 2.5,
which means that $G$ is $\sigma$-nilpotent. Hence  the first assertion of (4) is true. Finally, note that if, for example,  $P_{1}$ is not a maximal
subgroup of $P_{1}P_{2}$, then the chain $P_{1} < P_{1}P_{2}  < \cdots < P_{1} \cdots P_{n}=G$ can be refined to a maximal chain of $G$ of length $n$, at least. Hence $P_{1}$ is  $\sigma$-subnormal in $G$ by Lemma 4.5(ii). This contradiction shows that  $P_{1}$  forms an irreducible pair with $P_{i}$ for all $i > 1$.

(5)  {\sl The following assertions hold.}

(a) {\sl $P_{i}$ is elementary abelian for all $ i > 1$. Hence if $G$ possesses at least two non-$\sigma$-subnormal non-isomorphic Sylow subgroups, then all non-$\sigma$-subnormal Sylow subgroups are of prime order} (This follows from Lemma 6.1 and Claims (2) and (4)).

(b) {\sl If $P_{1}\leq H_{k}$ and $P_{1}$ is not of prime order, then $H_{1}, \ldots , H_{k-1}, H_{k+1}, \ldots , H_{t}$ are normal in $G$.}

Indeed, since $P_1$ is not $\sigma$-subnormal in $G$, $P_{1}$ is cyclic by Claim (2). Hence, if $i\ne k$ and  $P_{j}\leq H_{i}$, then $P_{j}$ does not form an irreducible pair with $P_{1} $ by Claim (a) and Lemma 6.1.  Therefore Claim (4) implies that  $P_{j}$ is $\sigma$-subnormal in $G$. This shows that every Sylow subgroup of $H_{i}$ is $\sigma$-subnormal in $G$. Hence $H_{i}$  is   normal in $G$ by Lemma 2.2(9)  and Proposition 2.5 for all $i\ne k$.

(c) {\sl If $P_{1}\leq  H_{k}$  and $V$ is the maximal subgroup  of $P_{1}$, then $|G:N_{G}(V)|$  is a $\sigma _{i_{k}}$-number. }

If $|P_{1}|$ is a prime, it is trivial.  Assume that  $|P_{1}|$ is not a prime and let $i\ne k$. Then $H_{i}$ is normal in $G$ by Claim (b).
On the other hand, Claim (2) implies that $V$ is $\sigma$-subnormal in $G$. Hence
$H_{i}\leq N_{G}(V)$ by Lemma 2.2(5). Hance we have (c).

(6) {\sl $D$ is a  Hall subgroup of $G$. Hence $D$ has a complement $M$ in $G$.}

Suppose that this is false and let $U$ be a Sylow  $p_{j}$-subgroup of $D$ such that $1 < U < P_{j}\leq H_{k}$.
Lemma 2.1 implies that  $$(G/N)^{{\mathfrak{N}}_{\sigma}}= G^{{\mathfrak{N}}_{\sigma}}N/N=DN/N$$
for any minimal normal subgroup $N$ of $G$.

Let $L$ be a minimal normal subgroup of $G$ contained in $D$.  Then $G/L$ is a group of type (ii). Indeed, assume that $G/L$ is $\sigma$-nilpotent and so $L=D$. Then $L < P_{j}$ and so $P_{1}$  does not form an irreducible pair with $P_{j}$. Hence $L < P_{1}=P_j$ by Claim (4).
From Claim (5)(b) it follows that $H_{1}\cdots H_{k-1}H_{k+1}\cdots H_{t}$ is normal in $G$, so $D\leq H_{1}\cdots H_{k-1}H_{k+1}\cdots H_{t}$ since
$G/H_{1}\cdots H_{k-1}H_{k+1}\cdots H_{t}\simeq H_{k}$ is $\sigma$-nilpotent, a contradiction. Hence $G/L$ is not $\sigma$-nilpotent, and so
$G/L$ is a group of type (ii) by Claim (1). Hence  $D/L$ is a Hall subgroup of $G/L$.
If $UL/L\ne 1$, then $UL/L$ is a Sylow  $p_{j}$-subgroup of $D/L$ and so $UL/L=P_{j}L/L$. Hence $P_{j}\leq UL\leq D$ and so $U = P_{j}$.
This contradiction shows that $UL/L=1$, so $U=L$.   Therefore $L < P_{j}$. But then, as above, we get that for any $i\ne k$ the subgroup $H_{i}$ is normal in $G$. Let $N$ be a minimal normal subgroup of $G$ contained in $H_{i}$. Then $G/N$ is not  $\sigma$-nilpotent since $i\neq k$, and so
$DN/N\simeq D$  is a Hall subgroup of $G/N$ by Claim (1), which implies that $U=P_{j}$ since $p_{j}\not \in \pi (N) \subseteq \pi (H_{i})$. This contradiction completes the proof of (6).

(7) {\sl Some Sylow subgroup $P_{i}^{x}$ of $G$ contained in $M$ is not $\sigma$-subnormal in $G$.}

Suppose that every Sylow subgroup $P_{i}^{x}$ of $G$ contained in $M$ is  $\sigma$-subnormal in $G$. Then $M$ is  $\sigma$-subnormal in $G$ by Proposition 2.5. Hence  there is a subgroup chain  $$M=M_{0} <  M_{1} < \cdots < M_{r}=G$$  such that  either $M_{i-1}\trianglelefteq M_{i}$ or
$M_{i}/(M_{i-1})_{M_{i}}$ is  ${\sigma}$-primary for all $i=1, \ldots , r$. Since $G/D=G/G^{\mathfrak{N}_{\sigma}}\simeq M$ is $\sigma$-nilpotent and $G$ is not $\sigma$-nilpotent, $M\ne G$. Hence we can assume without loss of generality that $M_{t-1} < G.$
 If $M_{t-1}$ is normal in $G$, then there is a  normal maximal subgroup $T$ of $G$ containing $M_{t-1}$ such that $G/T$ is nilpotent since $G$ is soluble. Then $D\leq T$, and so $G=DM=DT=T<G$, a contradiction.
Therefore $G/(M_{t-1})_{G}$ is $\sigma$-primary and thereby it is $\sigma$-nilpotent. But then $D\leq (M_{t-1})_G$, so $G=MD \leq M_{t-1} < G$. This contradiction shows  $M$ is not $\sigma$-subnormal in $G$. Hence we have (7).

(8) {\sl $D$ is  nilpotent}.

Suppose that $D$ is not nilpotent. Assume that  $G$ has a minimal normal subgroup $N\ne R$. Then in view of Claim (1), $D\simeq D/1=D/(R\cap N)$ is nilpotent.
Therefore $R$ is the unique minimal normal subgroup of $G$ and $R\nleq \Phi (G)$ by Lemma 2.1. It follows that  $R=C_{G}(R)\leq D$.

Since $D/R$ is nilpotent by Claim (1), there is a  normal subgroup $E/R$ of $D/R$ such that $R\leq E <  D$ and $D=E\rtimes P_{i}$ for some $i$.
 Assume that $P_{i}$ is $\sigma _{i_{k}}$-group.
The Frattini argument implies that for some $x\in G$ we have $M^{x}\leq N_{G}(P_{i})$.

Suppose that for each $r \ne k$  and for each Sylow subgroup $P$ of $G$ contained in $M^{x}$, where $P$ is  a $\sigma _{i_{r}}$-group,  we have $[P, P_{i}]=1$.
Then $G/E\simeq P_{i}M^{x}$ is $\sigma$-nilpotent and so $D\leq E$, a contradiction. Hence $G$ has a Sylow subgroup $P_{j}^{y}$ satisfying the following conditions:  $[P_{i}, P_{j}^{y}]\ne 1$,  $P_{j}^{y}\leq  M^{x}$ and $P_{j}^{y}$ is a $\sigma _{i_{r}}$-group for some $r \ne k$.

Then $P_{j}^{y}P_{i}$ is  $\sigma$-fiber and so $P_{j}^{y}$ is not $\sigma$-subnormal by Lemma 2.2(5) since $P_{j}^{y}\leq N_{G}(P_{i})$  and $[P_{i}, P_{j}^{y}]\ne 1$.
Now we show that  $P_{j}^{y}$ is of prime order. Assume that this is false. Then Claim (3) implies that $R$ is a $\sigma _{i_{r}}$-group, so $G$ has a minimal normal subgroup $N\ne R$ by Claim (5)(b), a contradiction. Hence $P_{j}^{y}$ is of prime order.

Note that since $P_j^y$ is not $\sigma$-subnormal, $P_{j}^{y}$ forms an irreducible pair with  $R$ and $P_{i}$ by Claim (4).  Then $C_{R}(P_{j}^{y})=1$ since $C_G(R)=R.$  Note also that  $C_{P_{i}}(P_{j}^{y})=1$. Indeed, since $P_{j}^{y}$ forms an irreducible pair with  $P_{i}$ and $P_{j}^{y}\leq C_{G}(P_{j}^{y})$, we have either $C_{P_{i}}(P_{j}^{y})=P_{i}$ or $C_{P_{i}}(P_{j}^{y})=1$. But the former  case is impossible since $[P_{i}, P_{j}^{y}]\ne 1$, so $C_{P_{i}}(P_{j}^{y})=1$.

Let $C=C_{G}(P_{j}^{y})\cap  RP_{i}$. Suppose that $C\ne 1$. Then $C= (R\cap C)U=U\leq  P_{i}^{a}$ for some $a \in R$, so
 $U\times  P_{j}^{y}\leq (P_{i}P_{j}^{y})^{b}\simeq (P_{i}P_{j}^{y})^{b}$ for some $b \in R$   since $G$ is soluble. It follows that
 $C=1$. Consequently, $P_j^y\cap C_G(RP_i)=1$. Hence $RP_{i}$ is nilpotent by  the Thompson theorem \cite[V, 8.14]{hupp}.  Thus  $P_{i}\leq C_{G}(R)=R$, a contradiction.  Therefore we have (8).

Claims (2)--(8) show that the necessity is true.

{\sl Sufficiency.}  If $G$ is of the type (i), it is clear. Now let $G$ be a group of the type (ii). Then $G$ is not $\sigma$-nilpotent by (ii)(b) and Proposition 3.4.
Hence $|\pi (G)|\leq m_{\sigma}(G)$ by Theorem 1.2(iii). Let $n=|\pi (G)|.$

In order ro prove that $m_{\sigma}(G)\leq |\pi (G)|$, we only need to prove that that every $n$-maximal subgroup of $G$ is $\sigma$-subnormal in $G$.
  Assume  that this is false and let $E$ be an $n$-maximal subgroup of $G$ such that $E$ is not $\sigma$-subnormal in $G$.
 Then  some  Sylow subgroup $E_{1}$ of $E$ is not $\sigma$-subnormal in $G$ by Proposition 2.5.
  We can assume that without loss of generality that
 $E_{1}\leq P_{1}$. Then $P_{1}$ is not $\sigma$-subnormal in $G$ by Lemma 2.2(4). If  $i > 1$  and $P_{1}P^{x}_{i}=P^{x}_{i}P_{1}$, then $P_{1} $ and $P^{x}_{i}$ are members of some Sylow basis of $G$ \cite[I, 4.16]{DH}.  By the hypothesis, $P_{1}$ forms an irreducible pair with $P^{x}_{i}$, $P_{1}$ is cyclic and the maximal subgroup of $P_{1}$ is $\sigma$-subnormal in $G$. Hence $E_{1}=P_{1}$, so every subgroup of $G$ containing $E_{1}$
 is a Hall subgroup of $G$ by Lemma 6.3. Then $E$ is exactly a $k$-maximal subgroup of $G$, where $k=|\pi (|G:E|)|$, by Lemma 6.4. Hence $k=n=|\pi (G)|
 $. But then $E=1$, so $E$ is $\sigma$-subnormal in $G$.  This contradiction completes the proof of the sufficiency.

 The theorem is proved.

\section{Final remarks and some open questions}

1. We say that  $G$ is  a group of \emph{$\sigma$-Spencer height $h_{\sigma}(G)=n$} if every maximal chain of $G$ of length $n$ contains a proper $\sigma$-subnormal entry and there exists at least one maximal chain of $G$ of length $n-1 $ which contains no any proper $\sigma$-subnormal entry.

In particular, if $\sigma =\{\{2\}, \{3\}, \cdots \}$,  we write $h(G)$ instead of $h_{\sigma}(G).$

It is clear that $h_{\sigma}(G)\leq  m_{\sigma}(G)$, and in general $m_{\sigma}(G)\ne  h_{\sigma}(G)$.

Theorem 1.2(ii)(iii) can be improved.

{\bf Theorem  7.1} \cite{prepr}.  {\sl Suppose that $G$ is $\sigma$-soluble Then the following statements hold.}

(i) {\sl   $l_{\sigma}(G)\leq h_{\sigma}(G)$.}

(ii) {\sl  If a soluble group $G$ is not $\sigma$-nilpotent, then $|\pi (G)|\leq h_{\sigma}(G)$.}

{\bf Corollary 7.2} (Spencer \cite{spencer}). {\sl Suppose that $G$ is a   soluble group.  Then:}

(i) {\sl $l(G) \leq h(G)$, where $l(G)$ is the nilpotent length of $G$.}

(ii) {\sl If  $h(G) <  |\pi (G)|$, then $G$ is nilpotent.}

2. In view of Theorems 1.2,  1.10  and 7.1, the following questions seam natural.

{\bf Question 7.1.} {\sl What is the structure of a soluble group  $G$ provided $|\pi (G)|=m_{\sigma}(G)+1$?}

{\bf Question 7.2.} {\sl What is the structure of a soluble group  $G$ provided $|\pi (G)|=h_{\sigma}(G)$?}

{\bf Question 7.3.} {\sl What is the structure of a soluble group  $G$ provided $|\pi (G)|=h_{\sigma}(G)+1$?}

Note that in the case when $\sigma =\{\{2\}, 3, \ldots \}$ the complete answers to these three  questions are known \cite{mann, spencer}.

\end{document}